\documentclass[]{article}
\begin{document}
\newtheorem{proposition}{Proposition}[section]
\newtheorem{definition}{Definition}[section]
\newtheorem{lemma}{Lemma}[section]

\title{\bf Sheaf Structures\\On a Class of Noncommutative Spectra}
\author{Keqin Liu\\Department of Mathematics\\The University of British Columbia\\Vancouver, BC\\
Canada, V6T 1Z2}
\date{September, 2011}
\maketitle

\begin{abstract} We introduce a class of noncommutative spectra and give the sheaf structure on the class of noncommutative spectra.
\end{abstract}

How to select a class of noncommutative rings such that we can rewrite algebraic geometry in the context of the class of  noncommutative rings? The purpose of this paper is to present our answer to this interesting problem. It is our belief that if a ring $R$ is in the right class of non-commutative rings which can be used to rewrite algebraic geometry over commutative rings, then $R$ should have  the following 
\begin{description}
\item[Partially Commutative Property:] $R$ is a graded ring, and for all homogeneous elements $x$ and $y$ of $R$, there exist homogeneous elements $x'$ and $y'$ of $R$ such that either $xy=yx'$ or $xy=y'x$.
\end{description}

In this paper, we introduce a class of noncommutative rings which are called Hu-Liu trirings by us. The trivial extension of a ring by a bimodule over the ring has been used in different mathematical areas for a long time (\cite{FGR} and \cite{HM}). Roughly speaking, a Hu-Liu triring is a trivial extension of a ring by a bimodule such that the bimodule carries an extra commutative ring structure.
Many results in commutative algebra and algebraic geometry have the satisfactory counterparts in Hu-Liu trirings.
The main result of this paper is the sheaf structures on the class of noncommutative spectra based on Hu-Liu trirings. In section 1, we give the basic properties of Hu-Liu trirings. The most important example of Hu-Liu trirings is triquaternions which is defined in section 1. Triquaternions, which are regarded as a kind of new numbers by us,  can be used to replace complex numbers to develop the counterpart of complex algebraic geometry. In section 2, we introduce prime triideals and prove some basic facts about prime triideals. In section 3, we use prime triideals to characterize the trinilradical. In section 4, we make the class of noncommutative spectra into a topological space by introducing extended Zariski topology. In section 5, we define localization of Hu-Liu trirings. In the last section of this paper, we explain how to define the  sheaf structures on the class of noncommutative spectra.

\medskip
Throughout this paper, the word ``ring'' means an associative ring with an identity. A ring $R$ is also denoted by $(R,\, +,\, \cdot)$ to indicate that $+$ is the addition and $\cdot$ is the multiplication in the ring $R$.

\bigskip
\section{Basic Definitions}

Let $A$ and $B$ be two subsets of a ring $(\, R, \, + , \, \cdot\,)$. We shall use $A+B$ and $AB$ to denote the following subsets of $R$ 
$$A+B:=\{\, a+b \, | \, a\in A, \, b\in B \, \}, \quad 
AB:=\{\, ab \, |\, a\in A, \, b\in B \,\}.$$

\medskip
We now introduce a class of noncommutative rings in the following 

\medskip
\begin{definition}\label{def2.1} A ring $R$ with a multiplication $\cdot$ is called a {\bf Hu-Liu triring} if the following three properties hold.
\begin{description}
\item[(i)] There exist a commutative subring $R_0$ of the ring $(R,\, +,\, \cdot)$ and a subgroups $R_1$ of the additive group $(R,\, +)$, called the {\bf even part} and {\bf odd part} of $R$ respectively, such that 
$R=R_0\oplus R_1$ (as Abelian groups) and 
\begin{equation}\label{eq1.1}
R_0R_0\subseteq R_0,\quad R_0R_1+R_1R_0\subseteq R_1, \quad R_1R_1=0;
\end{equation}
\item[(ii)] There exists a  binary operation $\sharp$ on the odd part $R_1$ such that 
$(\,R_1 , \, + , \, \,\sharp\, \,)$ is a commutative ring and the two associative products $\cdot$ and $\,\sharp\,$ satisfy the {\bf triassociative law}:
\begin{eqnarray}
\label{eq2} x (\alpha \,\sharp\, \beta) &=&(x \alpha) \,\sharp\, \beta, \\
\label{eq3} (\alpha \,\sharp\, \beta ) x &=&\alpha \,\sharp\, (\beta x), 
\end{eqnarray}
where $x\in R$ and $\alpha$, $\beta\in R_1$;
\item[(iii)] For each $x_0\in R_0$, we have
\begin{equation}\label{eq6}
R_1x_0=x_0R_1.
\end{equation}
\end{description}
\end{definition}

\medskip
A Hu-Liu triring $R=R_0\oplus R_1$, which clearly has the partially commutative property, is sometimes denoted by 
$(\,R=R_0\oplus R_1 , \, + , \, \cdot ,\,  \,\sharp\, \,)$,
where the associative product $\sharp$ on the odd part $R_1$ is called the {\bf local product}, and the identity $1^\sharp$ of the ring $(\, R_1 , \, + , \, \,\sharp\,\,)$ is called the {\bf local identity} of the triring $R$. If 
$(\,R=R_0\oplus R_1 , \, + , \, \cdot ,\,  \,\sharp\, \,)$ is a Hu-Liu triring, then $R$ is 
the trivial extension of a ring $R_0$ by a $R_0$-bimodule $R_1$, and the triassociative law interweaves the noncommutative ring structure on the  trivial extension with the commutative ring 
structure on the odd part.

\medskip
Since a commutative ring is a Hu-Liu triring with zero odd part, the concept of Hu-Liu trirings naturally generalizes the concept of commutative rings. The first and the most important example of Hu-Liu trirings which is not a commutative ring is triquaternions whose definition is given in the following example.

\medskip
\noindent
{\bf Example} Let 
${\bf Q}=\mathcal{R}\,1\oplus \mathcal{R}\,i\oplus\mathcal{R}\,j\oplus\mathcal{R}\,k$ be a 
$4$-dimensional real vector space, where $\mathcal{R}$ is the field of real numbers. Then 
$(\,{\bf Q}={\bf Q}_0\oplus {\bf Q}_1, \, +, \,\cdot, \, \,\sharp\, \,)$ is a Hu-liu triring, where ${\bf Q}_0=\mathcal{R}\,1\oplus \mathcal{R}\,i$ is the even part, 
${\bf Q}_1=\mathcal{R}\,j\oplus\mathcal{R}\,k$ is the odd part, the ring multiplication $\cdot$ and the local product $\sharp$ are defined by the following multiplication tables:
\begin{center}
{\huge\begin{tabular}{|c||c|c|c|c|}
\hline
$\cdot$&$1$&$i$&$j$&$k$\\
\hline\hline
$1$&$1$&$i$&$j$&$k$\\
\hline
$i$&$i$&$-1$&$k$&$-j$\\
\hline
$j$&$j$&$-k$&$0$&$0$\\
\hline
$k$&$k$&$j$&$0$&$0$\\
\hline
\end{tabular}}
\qquad\qquad
{\huge\begin{tabular}{|c||c|c|}
\hline
$\sharp$&$j$&$k$\\
\hline\hline
$j$&$j$&$k$\\
\hline
$k$&$k$&$-j$\\
\hline
\end{tabular}}
\end{center}
The Hu-Liu triring ${\bf Q}$ is called the {\bf triquaternions}.

\hfill\raisebox{1mm}{\framebox[2mm]{}}

\medskip
\begin{definition}\label{def2.2} Let $(\,R=R_0\oplus R_1, \, +, \,\cdot, \, \,\sharp\, \,)$ be a Hu-Liu triring with the identity $1$, and let $I$ be a subgroup of the additive group of $R$.
\begin{description}
\item[(i)] $I$ is called a {\bf triideal} of $R$ if $IR+RI\subseteq I$, 
$I=(R_0\cap I)\oplus (I\cap R_1)$, and $I\cap R_1$ is an ideal of the ring 
$(\,R_1, \, +, \, \,\sharp\, \,)$.
\item[(ii)] $I$ is called a {\bf subtriring} of $R$ if $1\in I$, $II\subseteq I$,  
$I=(R_0\cap I)\oplus (I\cap R_1)$ and $I\cap R_1$ is a subring of the ring 
$(\,R_1, \, +, \, \,\sharp\, \,)$.
\end{description}
\end{definition}

\medskip
Let $I$ be a triideal of a Hu-Liu triring $(\,R=R_0\oplus R_1, \, +, \,\cdot, \, \,\sharp\, \,)$. It is clear that $\displaystyle\frac{R}{I}=\left(\displaystyle\frac{R}{I}\right)_0\oplus 
\left(\displaystyle\frac{R}{I}\right)_1$  with
$\left(\displaystyle\frac{R}{I}\right)_i:=\displaystyle\frac{R_i+I}{I}$ for $i=0$ and $1$. We now define a local product on $\left(\displaystyle\frac{R}{I}\right)_1$ by
\begin{equation}\label{eq5}
(\alpha +I)\,\sharp\, (\beta +I):=\alpha\,\sharp\,\beta +I 
\quad\mbox{for $\alpha$, $\beta\in R_1$.}
\end{equation}
Then the local product defined by (\ref{eq5}) is well-defined, and the triassociative law holds. 
Therefore, $\displaystyle\frac{R}{I}$ becomes a Hu-Liu triring, which is called the {\bf quotient Hu-Liu triring} of $R$ with respect to the triideal $I$.

\medskip
\begin{definition}\label{def1.4} Let $R=R_0\oplus R_1$ and 
$\overline{R}=\overline{R}_0\oplus \overline{R}_1$ be Hu-Liu trirings. A map 
$\phi: R\to \overline{R}$ is called a {\bf triring homomorphism} if 
$$
\phi (x+y)=\phi (x)+\phi (y), \quad \phi (xy)=\phi (x)\phi (y), \quad
\phi (1_R)=1_{\overline{R}},
$$$$
\phi (R_0)\subseteq \overline{R}_0, \quad \phi (R_1)\subseteq \overline{R}_1,
$$$$
\phi (\alpha \,\sharp\, \beta )=\phi (\alpha) \,\sharp\, \phi (\beta), \quad
\phi(1^{\sharp})=\bar{1}^{\sharp}
$$
where $x$, $y\in R$, $\alpha$, $\beta\in R_1$, $1_R$ and $1_{\overline{R}}$ are the identities of $R$ and $\overline{R}$ respectively, and $1^{\sharp}$ and $\bar{1}^{\sharp}$ are the local identities of $R$ and $\overline{R}$ respectively. A bijective triring homomorphism is called a {\bf triring isomorphism}. We shall use $R\simeq \overline{R}$ to indicate that there exists a triring isomorphism from $R$ to $\overline{R}$.
\end{definition}

\medskip
Let  $\phi$ be a  triring homomorphism from a Hu-Liu triring $R$ to a Hu-Liu triring $\overline{R}$. the {\bf kernel} $Ker\phi$ and the {\bf image} $Im\phi$ of $\phi$ are defined by
$$ Ker\phi :=\{\, x \,|\, \mbox{$a\in R$ and $\phi (x)=0$} \, \}$$ 
and
$$ Im\phi :=\{\, \phi (x) \,|\, \mbox{$x\in R$} \, \}.$$ 
Clearly, $Ker\phi$ is a triideal of $R$, $Im\phi$ is a subtriring of $\overline{R}$ and 
$$\overline{\phi}: x+Ker\phi \to \phi (x) \quad\mbox{for $x\in R$}$$
is a triring isomorphism from the quotient Hu-Liu triring $\displaystyle\frac{R}{Ker\phi}$ to the subtriring $Im\phi$ of $\overline{R}$. If $I$ is triideal of $R$, then 
$$ \nu : x\mapsto x+I \quad\mbox{for $x\in R$}$$
is a surjective triring homomorphism from $R$ the quotient tiring $\displaystyle\frac{R}{I}$ with kernel $I$. The map $\nu$ is called the {\bf natural triring homomorphism}.

\bigskip
\begin{proposition}\label{pr1.1}  Let $\phi$ be a surjective homomorphism from a Hu-Liu triring $R=R_0\oplus R_1$ to a Hu-Liu triring $\overline{R}=\overline{R}_0\oplus \overline{R}_1$.
\begin{description}
\item[(i)] $\phi (R_i)=\overline{R}_i$ for $i=0$ and $1$.
\item[(ii)] Let
$$
\mathcal{S}: =\{\, I \,|\, \mbox{$I$ is a triideal of $R$ and $I \supseteq Ker \phi $} \, \}
$$
and
$$
\overline{\mathcal{S}}: =\{\, \overline{I} \,|\, \mbox{$\overline{I}$ is a triideal of $\overline{R}$} \, \}.
$$
The map 
$$
\Psi : I \mapsto \phi (I): =\{\, \phi (x) \,|\, x\in I \, \}
$$
is a bijection from $\mathcal{S}$ to $\overline{\mathcal{S}}$, and the inverse map $\Psi ^{-1} : \overline{\mathcal{S}} \to \mathcal{S}$ is given by
$$
\Psi ^{-1} : \overline{I} \mapsto \phi ^{-1} (\overline{I}) : = \{\, x \,|\, \mbox{$x\in R$ and $\phi (x)\in \overline{I}$ \, \}.}
$$
\item[(iii)] If $I$ is a triideal of $R$ containing $Ker \phi$, then the map
$$
x+I \mapsto \phi (x) + \phi (I) \quad\mbox{for $x\in R$}
$$
is a triring isomorphism from the quotient triring $\displaystyle\frac{R}{I}$ onto the quotient triring $\displaystyle\frac{\overline{R}}{\phi (I)}$.
\end{description}
\end{proposition}

\medskip
\noindent
{\bf Proof} A routine check.

\hfill\raisebox{1mm}{\framebox[2mm]{}}

\medskip
We now prove a basic property for Hu-Liu trirings.

\begin{proposition}\label{pr2.1}  Let $(\,R=R_0\oplus R_1, \, +, \,\cdot, \, \,\sharp\, \,)$ be a Hu-Liu triring. If $x_i\in R_i$  for $i\in\{\, 0, \, 1\,\}$, then both $Rx_0=R_0x_0\oplus R_1x_0$ and $R_1\,\sharp\, x_1$ are triideals of $R$.
\end{proposition}

\medskip
\noindent
{\bf Proof} Since $(R_1, \, +, \, \sharp)$ is a commutative ring, $R_1\,\sharp\, x_1$ is clearly a triideal of $R$.

\medskip
Using the triassociative law and (\ref{eq6}), we have
\begin{equation}\label{eq7}
\Big(R_i(R_0x_0+R_1x_0)\Big)\bigcup
\Big(R_0x_0+R_1x_0)R_i\Big)\subseteq R_0x_0+R_1x_0
\end{equation}
for $i=0$, $1$. By (\ref{eq7}), $R_0x_0\oplus R_1x_0$ is an ideal of the ring $(\,R, \, +, \, \cdot\,)$. Also, $R_1x_0=R_1\,\sharp\, (1^{\sharp}x_0)$ is obviously an ideal of the commutative ring $(\,R_1, \, +, \, \sharp\,)$. Thus, $R_0x_0\oplus R_1x_0$ is a triideal of $R$.

\hfill\raisebox{1mm}{\framebox[2mm]{}}

\medskip
The next position gives some operations about triideals in a Hu-Liu triring.

\medskip
\begin{proposition}\label{pr2.2}  Let $I$, $J$, $I_{\lambda}$ with $\lambda\in \Lambda$ be triideals of a Hu-Liu triring $R$.
\begin{description}
\item[(i)] The intersection $I\cap J$ and the sum 
$\displaystyle\sum_{\lambda\in \Lambda}I_{\lambda}$ are triideals of $R$.Moreover, we have
$(I\cap J)_i=I_i\cap J_i$ and 
$\left(\displaystyle\sum_{\lambda\in \Lambda}I_{\lambda}\right)_i=
\displaystyle\sum_{\lambda\in \Lambda}(I_{\lambda})_i$ for $i=0$, $1$.
\item[(ii)] The {\bf mixed product} 
$I\,\stackrel{\cdot}{\sharp}\,J:=(I\,\stackrel{\cdot}{\sharp}\,J)_0\oplus 
(I\,\stackrel{\cdot}{\sharp}\,J)_1$ of $I$ and $J$ is a triideal, where 
$(I\,\stackrel{\cdot}{\sharp}\,J)_0=I_0J_0$ and $(I\,\stackrel{\cdot}{\sharp}\,J)_1=I_1\,\sharp\,J_1$.
\end{description}
\end{proposition}

\medskip
\noindent
{\bf Proof} (i) It is clear.

\medskip
(ii) By triassociative law, we have
$$
R_0(I\,\stackrel{\cdot}{\sharp}\,J)\subseteq 
R_0I_0J_0+ R_0(I_1\,\sharp\,J_1)\subseteq I_0J_0+ (R_0I_1)\,\sharp\,J_1\subseteq I\,\stackrel{\cdot}{\sharp}\,J,
$$
$$
(I\,\stackrel{\cdot}{\sharp}\,J)R_0\subseteq 
I_0J_0R_0+(I_1\,\sharp\,J_1)R_0\subseteq I_0J_0+ I_1\,\sharp\,(J_1R_0)\subseteq I\,\stackrel{\cdot}{\sharp}\,J,
$$
$$
R_1(I\,\stackrel{\cdot}{\sharp}\,J)+(I\,\stackrel{\cdot}{\sharp}\,J)R_1\subseteq
R_1I_0J_0+I_0J_0R_1\subseteq I_1J_0+I_0J_1\subseteq I\,\stackrel{\cdot}{\sharp}\,J
$$
and $R_1\,\sharp\,(I_1\,\sharp\,J_1)\subseteq
R_1\,\sharp\,I_1\,\sharp\,J_1\subseteq I\,\stackrel{\cdot}{\sharp}\,J$. This proves that (ii) holds.

\hfill\raisebox{1mm}{\framebox[2mm]{}}

\bigskip
\section{Prime Triideals}

\medskip
We begin this section by introducing the notion of prime triideals.

\begin{definition}\label{def5.1} Let $(\, R=R_0\oplus R_1, \, + , \, \cdot\, , \,\sharp\, \,)$ be a Hu-Liu triring. An triideal $P=P_0\oplus P_1$ of $R$ is called a {\bf prime triideal} if $P\ne R$ and 
\begin{eqnarray}
\label{eq15}x_0y_0\in P_0  &\Rightarrow&  \mbox{$x_0\in P_0$ or $y_0\in P_0$},\\
\label{eq16}x_0y_1\in P_1  &\Rightarrow&  \mbox{$x_0\in P_0$ or $y_1\in P_1$},\\
\label{eq17}x_1y_0\in P_1  &\Rightarrow&  \mbox{$x_1\in P_1$ or $y_0\in P_0$},\\
\label{eq18}x_1\,\sharp\, y_1\in P_1 &\Rightarrow&  \mbox{$x_1\in P_1$ or $y_1\in P_1$},
\end{eqnarray}
where $x_i$, $y_i\in R_i$ for $i=0$ and $1$.
\end{definition}

\medskip
Let $(\, R=R_0\oplus R_1, \, + , \, \cdot\, \,\sharp\, \,)$ be a Hu-Liu triring. The set of all prime triideals of $R$ is called the {\bf trispectrum} of $R$ and denoted by $Spec^{\sharp}R$. It is clear that
$$Spec^{\sharp}R=Spec^{\sharp}_0R\cup Spec^{\sharp}_1R\quad\mbox{and}\quad 
Spec^{\sharp}_0R\cap Spec^{\sharp}_1R=\emptyset,$$
where
$$Spec^{\sharp}_0R:=\{\, P \, |\, \mbox{$P\in Spec^{\sharp}R$ and $P\supseteq R_1$}\,\}$$
is called the {\bf even trispectrum} of $R$ and 
$$Spec^{\sharp}_1R:=\{\, P \, |\, \mbox{$P\in Spec^{\sharp}R$ and $P\not\supseteq R_1$}\,\}$$
is called the {\bf odd trispectrum} of $R$.  

\medskip
Let $P=P_0\oplus P_1$ be a triideal of a Hu-Liu triring $R$. Then 
$P\in Spec^{\sharp}_0R$ if and only if $P_1=R_1$ and $P_0$ is a prime ideal of the commutative ring 
$(\, R_0, \, + , \, \cdot\, \,)$. It is also obvious that $P\in Spec^{\sharp}_1R$ if and only if
$P_0$ is a prime ideal of the commutative ring $(\, R_0, \, + , \, \cdot\, \,)$, $P_1$ is a prime ideal of the commutative ring $(\, R_1 , \, + , \, \,\sharp\, \, )$, and the $\left(\displaystyle\frac{R_0}{P_0}, \displaystyle\frac{R_0}{P_0}\right)$-bimodule $\displaystyle\frac{R}{P}$ is faithful as both left module and right module, where the left $\displaystyle\frac{R_0}{P_0}$-module action on $\displaystyle\frac{R}{P}$ is defined by 
$$
(x_0+P_0)(y+P):= x_0y+P \quad\mbox{for $x_0\in P_0$ and $y\in R$}
$$
and the right $\displaystyle\frac{R_0}{P_0}$-module action  on $\displaystyle\frac{R}{P}$ is defined by 
$$
(y+P)(x_0+P_0):= yx_0+P \quad\mbox{for $x_0\in P_0$ and $y\in R$}.
$$

\medskip
Clearly, the even trispectrum $Spec^{\sharp}_0R$ of a Hu-Liu triring $R$ is not empty. A basic property of Hu-Liu trirings is that the odd trispectrum $Spec^{\sharp}_1R$ is always not empty provided $R_1\ne 0$. This basic fact is a corollary of the following

\begin{proposition}\label{pr5.1} Let $(\, R=R_0\oplus R_1, \, + , \, \cdot\, \,\sharp\, \,)$ be a Hu-Liu triring with $R_1\ne 0$. If $P_1$ is a prime ideal of the commutative ring 
$(\, R_1 , \, + , \, \,\sharp\, \, )$, then there exists prime ideal $P_0$ of the commutative ring $(\, R_0 , \, + , \, \,\cdot\, \, )$ such that $P:=P_0\oplus P_1$ is a prime triideal of $R$, and $P_0$ contains every ideal $I_0$ of the ring  $(\, R_0 , \, + , \, \,\cdot\, \, )$ which has the property: $R_1I_0\subseteq P_1$.
\end{proposition}

\medskip
\noindent
{\bf Proof} Consider the set $\Omega$ defined by
$$
\Omega :=\{\,I_0\,|\, \mbox{$I_0$ is an ideal of $(\, R_0 , \, + , \, \,\cdot\, \, )$ and $R_1I_0\subseteq P_1$}\,\}.
$$

Clearly, $1\not\in I_0$ if $I_0\in \Omega$. Since $0\in \Omega$, $\Omega$ is not empty. The relation of inclusion, $\subseteq$, is a partial order on $\Omega$. Let $\Delta$ be a non-empty totally ordered subset of $\Omega$. Let $J_0:=\displaystyle\bigcup _{I_0\in\Delta} I_0$. Then $J_0\in\Omega$. Thus $J_0$ is an upper bound for $\Delta$ in $\Omega$. By Zorn's Lemma, the partial order set $(\Omega,\,\subseteq)$ has a maximal element $P_0$. We are going to prove that $P:=P_0\oplus P_1$ is a prime triideal of $R$. 
Clearly, $P=P_0\oplus P_1$ is a triideal satisfying (\ref{eq18}). Let $x_i$, $y_i\in R_i$ for $i=0$, $1$. 

\medskip
If $x_1y_0\in P_1$ and $x_1\not\in P_1$, then 
$x_1\,\sharp\,(1^{\sharp}y_0)=x_1y_0\in P_1$, which implies that $1^{\sharp}y_0\in P_1$. Hence, we have
\begin{eqnarray}\label{eq19}
R_1(P_0+R_0y_0)&\subseteq& R_1P_0+R_1R_0y_0\subseteq P_1+R_1y_0\nonumber\\
&=&P_1+R_1\,\sharp\,(1^{\sharp}y_0)\subseteq P_1+R_1\,\sharp\,P_1\subseteq P_1.
\end{eqnarray}
Using (\ref{eq19}) and the fact that $P_0+R_0y_0$ is an ideal of $R_0$, we get
$P_0+R_0y_0\in \Omega$. Since $P_0+R_0y_0\supseteq P_0$, we have to have 
$P_0+R_0y_0= P_0$, which implies that $y_0\in P_0$. This proves that (\ref{eq17}) holds.

\medskip
If $x_0y_1\in P_1$ and $y_1\not\in P_1$, then 
$(x_01^{\sharp})\,\sharp\,y_1=x_0y_1\in P_1$, which implies that $x_01^{\sharp}\in P_1$. It follows from this fact and (\ref{eq6}) that
\begin{eqnarray}\label{eq20}
R_1(P_0+R_0x_0)&\subseteq& R_1P_0+R_1R_0x_0\subseteq P_1+R_1x_0= P_1+x_0R_1\nonumber\\
&=&P_1+(x_01^{\sharp})\,\sharp\,R_1\subseteq P_1+P_1\,\sharp\,R_1\subseteq P_1.
\end{eqnarray}
Using (\ref{eq20}) and the fact that $P_0+R_0x_0$ is an ideal of $R_0$, we get
$P_0+R_0x_0\in \Omega$. Since $P_0+R_0x_0\supseteq P_0$, we have to have 
$P_0+R_0x_0= P_0$, which implies that $x_0\in P_0$. This proves that (\ref{eq16}) holds.

\medskip
Finally, if $x_0y_0\in P_0$, then 
$(1^{\sharp}x_0)\,\sharp\,(1^{\sharp}y_0)=1^{\sharp}x_0y_0\in P_1$. Thus, 
$(1^{\sharp}x_0)\in P_1$ or $(1^{\sharp}y_0)\in P_1$, which implies that either $x_0\in P_0$ or 
$y_0\in P_0$ by the argument above. This proves that (\ref{eq15}) holds.

\medskip
Summarizing what we have proved, $P=P_0\oplus P_1$ is a prime triideal of $R$. 

\hfill\raisebox{1mm}{\framebox[2mm]{}}

\medskip
The following proposition gives another characterization of prime triideals.

\begin{proposition}\label{pr5.2} Let $(\, R=R_0\oplus R_1, \, + , \, \cdot , \, \,\sharp\, \,)$ be a Hu-Liu triring. The following are equivalent.
\begin{description}
\item[(i)] $P$ is a prime triideal.
\item[(ii)] For two triideals $I$, $J$ of $R$, $I\,\stackrel{\cdot}{\sharp}\,J\subseteq P$ implies that $I\subseteq P$ or $J\subseteq P$.
\end{description}
\end{proposition}

\medskip
\noindent
{\bf Proof} (i) $\Rightarrow$ (ii): Assume that $I\,\stackrel{\cdot}{\sharp}\,J\subseteq P$ and $I\not\subseteq P$. Then either $I_0\not\subseteq P_0$ or $I_1\not\subseteq P_1$. Let 
$y_i\in J_i$ with $i=0$, $1$.

\medskip
If $I_0\not\subseteq P_0$, then there exists $x_0\in P_0\setminus I_0$. Clearly, 
$x_0y_i\in I\,\stackrel{\cdot}{\sharp}\,J\subseteq P$ for $i=0$, $1$. Since $P$ is a prime triideal, we get $y_i\in P$ for $i=0$, $1$ by (\ref{eq15}) and (\ref{eq16}). Thus $J\subseteq P$ in this case. 

\medskip
If $I_1\not\subseteq P_1$, then there exists $x_1\in P_1\setminus I_1$. Since 
$x_1y_0\in I\,\stackrel{\cdot}{\sharp}\,J\subseteq P$ and 
$x_1\,\sharp\, y_1\in \,\stackrel{\cdot}{\sharp}\,J\subseteq P$, we get $y_i\in P$ for $i=0$, $1$ by (\ref{eq17}) and (\ref{eq18}). Thus we also get $J\subseteq P$ in this case. 

\bigskip
(ii) $\Rightarrow$ (i): Let $x_i$, $y_i\in P_i$ with for $i=0$, $1$. By Proposition \ref{pr2.1}, both $R_0x_0\oplus R_1x_0$ and $R_0y_0\oplus R_1y_0$ are triideals. By the definition of the mixed product of two triideals in Proposition \ref{pr2.2}, we have
$$(R_0x_0\oplus R_1x_0)\,\stackrel{\cdot}{\sharp}\,(R_0y_0\oplus R_1y_0)
=R_0x_0R_0y_0+(R_1x_0)\,\sharp\, (R_1y_0).$$ 
Since $(R_1x_0)\,\sharp\, (R_1y_0)=R_1x_0y_0$, we get
\begin{equation}\label{eq21}
(R_0x_0\oplus R_1x_0)\,\stackrel{\cdot}{\sharp}\,(R_0y_0\oplus R_1y_0)
\subseteq R_0x_0y_0+R_1x_0y_0.
\end{equation}

\medskip
If $x_0y_0\in P_0$ and $x_0\not\in P_0$, then 
\begin{equation}\label{eq22}
R_0x_0\oplus R_1x_0\not\subseteq P.
\end{equation}
It follows from  (\ref{eq21}) and (\ref{eq22}) that $R_0y_0\oplus R_1y_0\subseteq P$, which implies that $y_0\in P_0$. Thus, (\ref{eq15}) holds.

\medskip
If $x_0y_1\in P_1$ and $x_0\not\in P_0$, then (\ref{eq22}) holds.
By (\ref{eq6}), $1^{\sharp}x_0=x_0z_1$ for some $z_1\in R_1$. Hence, we have
\begin{eqnarray}\label{eq23}
&&(R_0x_0\oplus R_1x_0)\,\stackrel{\cdot}{\sharp}\,(R_1\,\sharp\, y_1)= (R_1x_0)\,\sharp\,(R_1\,\sharp\, y_1)\subseteq 
R_1\,\sharp\, (1^{\sharp}x_0)\,\sharp\, y_1\nonumber\\
&\subseteq &R_1\,\sharp\, (x_0z_1)\,\sharp\, y_1
\subseteq R_1\,\sharp\, (x_01^{\sharp})\,\sharp\,z_1\,\sharp\, y_1
\subseteq R_1\,\sharp\, z_1\,\sharp\,(x_01^{\sharp})\,\sharp\, y_1\nonumber\\
&\subseteq& R_1\,\sharp\, \big(x_0(1^{\sharp}\,\sharp\,y_1)\big)
\subseteq R_1\,\sharp\, (x_0y_1).
\end{eqnarray}
It follows from  (\ref{eq21}) and (\ref{eq23}) that $R_1\,\sharp\, y_1\subseteq P$, which implies that $y_1\in P$. Thus, (\ref{eq16}) holds.

\medskip
If $x_1y_0\in P_1$ and $x_1\not\in P_1$, then
\begin{equation}\label{eq24} 
R_1\,\sharp\, x_1\not\subseteq P.
\end{equation}

It follows from $x_1y_0\in P_1$ that
\begin{eqnarray}\label{eq25}
&&(R_1\,\sharp\, x_1)\,\stackrel{\cdot}{\sharp}\,(R_0y_0+R_1y_0)
=(R_1\,\sharp\, x_1)\,\sharp\,(R_1y_0)=R_1\,\sharp\, x_1\,\sharp\,R_1\,\sharp\,(1^{\sharp}y_0)
\nonumber\\
&=&R_1\,\sharp\, x_1\,\sharp\,(1^{\sharp}y_0)=R_1\,\sharp\, (x_1\,\sharp\,1^{\sharp})y_0
=R_1\,\sharp\, (x_1y_0)\subseteq P.
\end{eqnarray}
By (\ref{eq24}) and (\ref{eq25}), we get $R_0y_0+R_1y_0\subseteq P$, which implies that 
$y_0\in P$. Thus, (\ref{eq17}) holds.

\medskip
If $x_1\,\sharp\,y_1\in P_1$ and $x_1\not\in P_1$, then (\ref{eq24}) holds and
\begin{equation}\label{eq26} 
(R_1\,\sharp\, x_1)\,\stackrel{\cdot}{\sharp}\,(R_1\,\sharp\, y_1)
=R_1\,\sharp\, x_1\,\sharp\,y_1\subseteq P.
\end{equation}
By (\ref{eq24}) and (\ref{eq26}), $R_1\,\sharp\, y_1\subseteq P$, which implies that 
$y_1\in P$. Thus, (\ref{eq18}) holds.

\medskip
This proves that $P$ is a prime triideal.

\hfill\raisebox{1mm}{\framebox[2mm]{}}

\bigskip
\section{Trinilradicals}

Let $R=R_0\oplus R_1$ be a Hu-Liu triring with the local identity $1^\sharp$. For 
$\alpha \in R_1$, the {\bf local $n$th power} $\alpha ^{\sharp n}$ is defined by:
$$
\alpha^{\sharp n} :=\left\{ \begin{array}{ll}
1^\sharp, & \qquad\mbox{if $n=0$;}\\
\underbrace{\alpha \,\sharp\, \alpha  \,\sharp\, \cdots  \,\sharp\, \alpha}_n , & 
\qquad\mbox{if $n$ is a positive integer.} \end{array}\right.
$$
The products $(x^m)(\alpha^{\sharp n})$ and $(\alpha^{\sharp n})(x^m)$ will be denoted by $x^m\alpha^{\sharp n}$ and $\alpha^{\sharp n}x^m$ respectively, where $x\in R$ and $\alpha \in R_1$.

\begin{proposition}\label{pr6.1} Let $(\, R=R_0\oplus R_1, \, + , \, \cdot , \, \,\sharp\, \,)$ be a Hu-Liu triring with the local identity $1^\sharp$. 
If $x$, $y\in R$, $\alpha$, $\beta\in R_1$ and $m\in\mathcal{Z}_{>0}$, then
\begin{equation}\label{eq27}
(x\alpha)\,\sharp\, (y\beta)=(xy)(\alpha\,\sharp\, \beta), \qquad
(\alpha x)\,\sharp\, (\beta y)=(\alpha\,\sharp\, \beta)xy
\end{equation}
and
\begin{equation}\label{eq28}
(x\alpha)^{\sharp m}=x^m \alpha ^{\sharp m}, \qquad 
(\alpha x)^{\sharp m}=\alpha ^{\sharp m}x^m.
\end{equation}
\end{proposition}

\medskip
\noindent
{\bf Proof} By the triassociative law, we have
$$
(x\alpha )\,\sharp\, (y\beta)=x\left(\alpha\,\sharp\, (y\beta)\right)=
x\left((y\beta)\,\sharp\, \alpha\right)=(xy)(\beta\,\sharp\, \alpha)=(xy)(\alpha\,\sharp\, \beta)
$$
and
$$
(\alpha x)\,\sharp\, (\beta y)=\left((\alpha x)\,\sharp\, \beta)\right)y=
\left((\beta)\,\sharp\, (\alpha x)\right)y=(\beta\,\sharp\, \alpha)xy=
(\alpha\,\sharp\, \beta)xy.
$$
Hence, (\ref{eq27}) holds. Clearly, (\ref{eq28}) follows from (\ref{eq27}).

\hfill\raisebox{1mm}{\framebox[2mm]{}}

\medskip
\begin{definition}\label{def6.1} Let $(\, R=R_0\oplus R_1, \, + , \, \cdot , \, \,\sharp\, \,)$ be a Hu-Liu triring. 
\begin{description}
\item[(i)] An element $x$ of $R$ is said to be {\bf  trinilpotent} if 
\begin{equation}\label{eq29}
\mbox{$x_0^m=0$ and $x_1^{\,\sharp\, n}=0$ for some $m$, $n\in\mathcal{Z}_{>0}$,}
\end{equation}
where $x_0$ and $x_1$ are the even component and the old component of $x$ respectively.
\item[(ii)] The set of all graded nilpotent elements of $R$ is called the {\bf trinilradical} of $R$ and denoted by $nilrad ^{\sharp}(R)$ or $\sqrt[\sharp]{0}$.
\end{description}
\end{definition}

\medskip
The ordinary nilradical of a ring $(\, A, \, + , \, \cdot  \,)$ is denoted by $nilrad(A)$ or
$nilrad(A, + , \cdot)$; that is,
$$
nilrad(A):=\{\, x \,|\, \mbox{$x^m=0$ for some $m\in\mathcal{Z}_{>0}$}\,\}.
$$

If $(\, R=R_0\oplus R_1, \, + , \, \cdot , \, \,\sharp\, \,)$ is a Hu-Liu triring, then
$$
nilrad(R, +, \cdot)=nilrad(R_0, + , \cdot)\oplus R_1
$$
and
\begin{equation}\label{eq30}
nilrad^{\sharp}R=nilrad(R_0, + , \cdot)\oplus nilrad(R_1, +, \,\sharp\,).
\end{equation}
Hence, $nilrad^{\sharp}R\subseteq nilrad(R, +, \cdot)$, i.e., the trinilradical of a Hu-Liu triring $R$ is smaller than the ordinary nilradical of the ring $(R,\,+, \,\cdot)$.

\begin{proposition}\label{pr6.2} Let $(\, R=R_0\oplus R_1, \, + , \, \cdot , \, \,\sharp\, \,)$ be a Hu-Liu triring.
\begin{description}
\item[(i)] The trinilradical $nilrad ^{\sharp}(R)$ is a triideal of $R$.
\item[(ii)] $nilrad ^{\sharp}\left( \displaystyle\frac{R}{nilrad ^{\sharp}(R)}\right)=0.$
\end{description}
\end{proposition}

\medskip
\noindent
{\bf Proof} (i) By the definition of trinilradicals, we have 
$$
(nilrad^{\sharp}R)\cap R_0=nilrad(R_0, + , \cdot)\quad\mbox{and}\quad
(nilrad^{\sharp}R)\cap R_1=nilrad(R_1, +, \,\sharp\,).
$$
Using the fact above and (\ref{eq30}), we need only to prove
\begin{equation}\label{eq31}
\{xa,\, ax\}\subseteq nilrad ^{\sharp}R\quad\mbox{for $x\in R$ and $a\in nilrad ^{\sharp}R$}.
\end{equation}

Let $x=x_0+x_1\in R$ and $a=a_0+a_1\in nilrad ^{\sharp}R$, where $x_i$, $a_i\in R_i$ for $i=0$, $1$. Then we have $a_0^m=a_1^{\,\sharp\, m}=0$ for some $m\in\mathcal{Z}_{>0}$. Since
\begin{eqnarray}
\label{eq32} xa&=&(x_0+x_1)(a_0+a_1)=x_0a_0+(x_0a_1+x_1a_0),\\
(x_0a_0)^m&=&x_0^ma_0^m=x_0^m0=0,\nonumber\\
(x_0a_1)^{\sharp m}&=&x_0^ma_1^{\sharp m}=x_0^m0=0,\nonumber\\
(x_1a_0)^{\sharp m}&=&(x_1a_0)^{\sharp m}=x_1^{\sharp m}a_0^m=x_1^{\sharp m}0=0,\nonumber
\end{eqnarray}
we get
\begin{equation}\label{eq33}
x_0a_0\in nilrad(R_0, \,+, \,\cdot) \quad\mbox{and}\quad 
x_0a_1+x_1a_0\in nilrad ^{\sharp}(R_1, \,+,\, \sharp).
\end{equation}

It follows from (\ref{eq32}) and (\ref{eq33}) that $xa\in nilrad ^{\sharp}R$. Similarly, we have $ax\in nilrad ^{\sharp}R$. This proves (i).

\medskip
(ii) If $x+nilrad ^{\sharp}R\in \displaystyle\frac{R}{nilrad ^{\sharp}R}$, then there exist positive integers $m$ and $n$ such that
$$x_0^m+nilrad ^{\sharp}R=\left(x_0+nilrad ^{\sharp}R\right)^m=nilrad ^{\sharp}R$$
and
$$x_1^{\,\sharp\, n}+nilrad ^{\sharp}R=\left((x_1+nilrad ^{\sharp}R)_1\right)^{\,\sharp\, n}=nilrad ^{\sharp}R,$$
where $x=x_0+x_1$ and $x_i\in R_i$ for $i=0$, $1$. Hence, we get
$$ x_0^m\in nilrad ^{\sharp}R \quad\mbox{and}\quad x_1^{\,\sharp\, n}\in nilrad ^{\sharp}R,$$
which imply that
$$x_0^{mu}=(x_0^m)^u=0  \quad\mbox{and}\quad 
x_1^{\,\sharp\, (nv)}=(x_1^{\,\sharp\, n})^{\,\sharp\, v}=0$$
for some $u$, $v\in \mathcal{Z}_{>0}$. Thus, $x=x_0+x_1\in nilrad ^{\sharp}R$ or 
$x+nilrad ^{\sharp}R$ is the zero element of $\displaystyle\frac{R}{nilrad ^{\sharp}R}$. This proves (ii).

\hfill\raisebox{1mm}{\framebox[2mm]{}}

\medskip
If $I$ is a triideal of a Hu-Liu triring $R=R_0\oplus R_1$, then the {\bf trinilradical} $\sqrt[\sharp]{I}$ of $I$ is defined by
$$
\sqrt[\sharp]{I}:=\{\, x\in R \, | \, \mbox{$x_0^m\in I_0$ and $x_1^{\,\sharp\, n}\in I_1$ for some $m$, $n\in\mathcal{Z}_{>0}$} \,\},
$$
where $x=x_0+x_1$, $x_i\in R_i$ and $I_i=I\cap R_i$ for $i=0$, $1$. Since
$$
nilrad ^{\sharp}\left(\displaystyle\frac{R}{I}\right)=\displaystyle\frac{\sqrt[\sharp]{I}}{I},
$$
$\sqrt[\sharp]{I}$ is a triideal of $R$. A triideal $I$ of a Hu-Liu triring $R$ is called a {\bf radical triideal} if $\sqrt[\sharp]{I}=I$.

\medskip
We now characterize the trinilradical of a Hu-Liu triring by using prime triideals. 

\begin{proposition}\label{pr6.3} Let $(\, R=R_0\oplus R_1, \, + , \, \cdot\, , \,\sharp\, \,)$ be a Hu-Liu triring. The trinilradical of $R$ is the intersection of the prime triideals of $R$.
\end{proposition}

\medskip
\noindent
{\bf Proof}  Let $x=x_0+x_1$ be any element of $nilrad ^{\sharp}(R)$, where $x_0\in R_0$ and $x_1\in R_1$. Then $x_0^m=0$ and $x_1^{\,\sharp\, n}=0$ for some $m$, $n\in\mathcal{Z}_{>0}$. Let $P=P_0\oplus P_1$ be any prime triideal of $R$. Since $x_0^m=0\in P_0$, we have 
\begin{equation}\label{eq34} 
x_0\in P_0
\end{equation}
by (\ref{eq15}).

If $P\not\supseteq R_1$, then $P_1$ is a prime ideal of the commutative ring 
$(\, R_1, \, + , \, \,\sharp\, \,)$. Using this fact and $x_1^{\,\sharp\, n}=0$, we get
\begin{equation}\label{eq35} 
x_1\in P_1.
\end{equation}
If $P\supseteq R_1$, then (\ref{eq35}) is obviously true. It follows from (\ref{eq34}) and (\ref{eq35}) that $x=x_0+x_1\in P$. This proves that
\begin{equation}\label{eq36} 
nilrad ^{\sharp}(R)\subseteq \bigcap_{P\in Spec^{\sharp}R} P.
\end{equation}

\medskip
Conversely, we prove that
\begin{equation}\label{eq37} 
z\not\in nilrad ^{\sharp}(R) \,\Rightarrow\, z\not\in \bigcap_{P\in Spec^{\sharp}R} P.
\end{equation}

\medskip
\underline{\it Case 1:} $z^m\ne 0$ for all $m\in \mathcal{Z}_{>0}$, in which case, 
$z^m\not\in R_1$ for all $m\in \mathcal{Z}_{>0}$. Hence, $z+ R_1$ is not a nilpotent element of the commutative ring $\displaystyle\frac{R}{R_1}$; that is,
$$
z+ R_1\ne nilrad\left(\displaystyle\frac{R}{R_1} \right)
=\bigcap _{\frac{I}{R_1}\in Spec \left( \frac{R}{R_1}\right)} 
\left(\displaystyle\frac{I}{R_1}\right),
$$
where $Spec \left( \displaystyle\frac{R}{R_1}\right)$ is the ordinary spectrum of the commutative ring
$\left(\displaystyle\frac{R}{R_1},\, +,\, \cdot\right)$.
Hence, there exists a prime ideal $\displaystyle\frac{I}{R_1}$ of the commutative ring $\displaystyle\frac{R}{R_1}$ such that $z\not\in I$. Since $I$ is a prime triideal of $R$, (\ref{eq37}) holds in this case.

\medskip
\underline{\it Case 2:} $z^m= 0$ for some $m\in \mathcal{Z}_{>0}$. Let $z=z_0+z_1$ with 
$z_0\in R_0$ and $z_1\in R_1$. Then $0=z^m=z_0^m+z_0^{m-1}r_1$ for some $r_1\in R_1$ by (\ref{eq6}). Thus, $z_0^m=0$, which implies that  $z_1^{\,\sharp\, n}\ne 0$ for all 
$n\in \mathcal{Z}_{>0}$ in this case.
We  now consider the following set
$$
T:=\left\{\, J \, \left| \, \mbox{$J$ is a triideal of $R$ and $z_1^{\,\sharp\, n}\not\in J$ for all $n\in \mathcal{Z}_{>0}$}\right.\right\}.
$$
Since $\{0\}\in T$, $T$ is nonempty. Clearly, $(\,T, \, \subseteq \,)$ is a partially order set, where $\subseteq $ is the relation of set inclusion. If $\{\, J_{\lambda} \, |\, \lambda\in \Lambda \,\}$ is a nonempty totally ordered subset of $T$, then 
$\bigcup_{\lambda\in \Lambda} J_{\lambda}$ 
is an upper bound of 
$\{\, J_{\lambda} \, |\, \lambda\in \Lambda \,\}$ in $T$. By Zorn's Lemma, the partially ordered set $(\,T, \, \subseteq \,)$ has a maximal element $P$. We are going to prove that $P$ is a prime triideal of $R$.

\medskip
Let $x=x_0+x_1$ and $y=y_0+y_1$ be two elements of $R$, where $x_i$, $y_i\in R_i$ for $i=0$, $1$. First, if $x_0\not\in P_0$ and $y_0\not\in P_0$, then 
\begin{equation}\label{eq38}
P\subset P+Rx_0\quad\mbox{and}\quad P\subset P+Ry_0.
\end{equation}
Since both $P+Rx_0$ and $P+Ry_0$ are triideals of $R$ by Proposition~\ref{pr2.1}, (\ref{eq38}) implies that
$$
z_1^{\,\sharp\, u}\in P+Rx_0\quad\mbox{and}\quad z_1^{\,\sharp\, v}\in P+Ry_0
$$
or
$$
z_1^{\,\sharp\, u}\in P_1+R_1x_0\quad\mbox{and}\quad z_1^{\,\sharp\, v}\in P_1+R_1y_0\quad\mbox{for some $u$, $v\in \mathcal{Z}_{>0}$.}
$$
Thus, we have
\begin{eqnarray*}
&&z_1^{\,\sharp\, (u+v)}=z_1^{\,\sharp\, u}\,\sharp\, z_1^{\,\sharp\, v}\in (P_1+R_1x_0)\,\sharp\, (P_1+R_1y_0)\\&&\\
&\subseteq& 
\underbrace{P_1\,\sharp\, P_1+ P_1\,\sharp\, (R_1y_0)+
(R_1x_0)\,\sharp\, P_1}_{\mbox{This is a subset of $P$}}+ (R_1x_0)\,\sharp\, (R_1y_0)\\&&\\
&\subseteq& P+R_1x_0y_0,
\end{eqnarray*}
which implies that $x_0y_0\not\in P$. This proves that
\begin{equation}\label{eq39} 
\mbox{$x_0\not\in P_0$ and $y_0\not\in P_0$} \,\Rightarrow\, x_0y_0\not\in P_0.
\end{equation}

Next, if $x_0\not\in P_0$ and $y_1\not\in P_1$, then 
\begin{equation}\label{eq40}
P\subset P+Rx_0\quad\mbox{and}\quad P\subset P+R_1\,\sharp\,y_1.
\end{equation}
Since both $P+Rx_0$ and $P+R_1\,\sharp\,y_1$ are triideals of $R$ by Proposition~\ref{pr2.1}, (\ref{eq40}) implies that there exist some positive integers $s$ and $t$ such that
$$
z_1^{\,\sharp\, s}\in P+Rx_0\quad\mbox{and}\quad z_1^{\,\sharp\, t}\in P+ R_1 \,\sharp\, y_1
$$ 
or
\begin{equation}\label{eq41} 
z_1^{\,\sharp\, s}\in P_1+R_1x_0\quad\mbox{and}\quad z_1^{\,\sharp\, t}\in 
P_1+ R_1 \,\sharp\, y_1.
\end{equation}
It follows that
\begin{eqnarray*}
&&z_1^{\,\sharp\, (s+t)}=z_1^{\,\sharp\, s}\,\sharp\, z_1^{\,\sharp\, t}\in (P_1+R_1x_0)\,\sharp\, (P_1+ R_1 \,\sharp\, y_1)\nonumber\\&&\nonumber\\
&\subseteq& 
\underbrace{P_1\,\sharp\, P_1+ P_1\,\sharp\, R_1 \,\sharp\, y_1+
(R_1x_0)\,\sharp\, P_1}_{\mbox{This is a subset of $P_1$}}
+ (R_1x_0)\,\sharp\, R_1\,\sharp\, y_1\nonumber\\
&\subseteq& P_1+ R_1\,\sharp\,(x_0y_1),
\end{eqnarray*}
which implies that $x_0y_1\not\in P_1$. This proves that
\begin{equation}\label{eq44} 
\mbox{$x_0\not\in P_0$ and $y_1\not\in P_1$} \,\Rightarrow\, x_0y_1\not\in P_1.
\end{equation}

\medskip
Similarly, we have
\begin{equation}\label{eq47} 
\mbox{$y_1\not\in P_1$ and $x_0\not\in P_0$} \,\Rightarrow\, y_1x_0\not\in P_1
\end{equation}
and
\begin{equation}\label{eq48} 
\mbox{$x_1\not\in P_1$ and $y_1\not\in P_1$} \,\Rightarrow\, x_1y_1\not\in P_0.
\end{equation}

\medskip
By (\ref{eq39}), (\ref{eq44}), (\ref{eq47}) and (\ref{eq48}), $P$ is a prime triideal. Since $z_1\not\in P$, (\ref{eq37}) also holds in Case 2.

\medskip
It follows from (\ref{eq36}) and (\ref{eq37}) that Proposition~\ref{pr6.3} is true.

\hfill\raisebox{1mm}{\framebox[2mm]{}}

\bigskip
The next proposition is a corollary of Proposition~\ref{pr6.3}.

\begin{proposition}\label{pr6.4} If $I$ is a triideal of a Hu-Liu triring $R$ and $I\ne R$, then
$$\sqrt[\sharp]{I}=\bigcap _{\mbox{$P\in Spec^{\sharp}R$ and $P\supseteq I$}} P.$$
\end{proposition}

\medskip
\noindent
{\bf Proof} By Proposition~\ref{pr6.3}, we have
\begin{eqnarray*}
x\in \sqrt[\sharp]{I}
&\Leftrightarrow& x+I\in nilrad ^{\sharp}\left(\frac{R}{I}\right)
=\bigcap _{\mbox{$\frac{P}{I}\in Spec^{\sharp}\left(\frac{R}{I}\right)$}}\frac{P}{I}\\
&\Leftrightarrow& x\in \bigcap _{\mbox{$P\in Spec^{\sharp}R$ and $P\supseteq I$}} P.
\end{eqnarray*}

\hfill\raisebox{1mm}{\framebox[2mm]{}}

\bigskip
\section{Extended Zariski Topology}

Let $(\, R=R_0\oplus R_1, \, + , \, \cdot , \, \,\sharp\, \,)$ be a Hu-Liu triring. For a triideal $I$ of $R$, we define a subset $\mathcal{V}(I)$ of $Spec^{\sharp}R$ by
\begin{equation}\label{eq49} 
V^{\sharp}(I): =\{\, P \, |\, \mbox{$P\in Spec^{\sharp}R$ and $P\supseteq I$} \,\}.
\end{equation}

\begin{proposition}\label{pr7.1} Let $R$ be a Hu-Liu triring.
\begin{description}
\item[(i)] $V^{\sharp}(0)=spec^{\sharp}R$ and $V^{\sharp}(R)=\emptyset$.
\item[(ii)] $V^{\sharp}(I)\cup V^{\sharp}(J)
=V^{\sharp}(I\cap J)=V^{\sharp}(I\,\stackrel{\cdot}{\sharp}\,J)$, where $I$ and $J$ are two triideals of $R$.
\item[(iii)] $\displaystyle\bigcap _{\lambda \in \Lambda} 
V^{\sharp}\left(I_{(\lambda)}\right)=
V^{\sharp}\left(\displaystyle\sum _{\lambda \in \Lambda}I_{(\lambda)}\right)$, where 
$\left\{\left.\,I_{(\lambda)}\, \right| \, \lambda \in \Lambda\, \right\}$ is a set of triideals of $R$.
\end{description}
\end{proposition}

\medskip
\noindent
{\bf Proof} Since (i) and (iii) are clear, we need only to prove (ii).

\medskip
Since $I\,\stackrel{\cdot}{\sharp}\,J\subseteq I\cap J\subseteq I$, we get $V^{\sharp}(I\,\stackrel{\cdot}{\sharp}\,J)\supseteq V^{\sharp}(I\cap J)
\supseteq V^{\sharp}(I)$. Similarly, we have $V^{\sharp}(I\,\stackrel{\cdot}{\sharp}\,J)\supseteq V^{\sharp}(I\cap J)\supseteq V^{\sharp}(J)$. Thus, we get
\begin{equation}\label{eq50} 
V^{\sharp}(I)\cup V^{\sharp}(J)\subseteq V^{\sharp}(I\cap J)
\subseteq V^{\sharp}(I\,\stackrel{\cdot}{\sharp}\,J).
\end{equation}

\medskip
Conversely, if $P\in V^{\sharp}(I\,\stackrel{\cdot}{\sharp}\,J)$, then $I\,\stackrel{\cdot}{\sharp}\,J\subseteq P$. By Proposition~\ref{pr5.2}, we get $I\subseteq P$ or 
$J\subseteq P$. Hence, $P\in V^{\sharp}(I)\cup V^{\sharp}(J)$. This proves that
\begin{equation}\label{eq51} 
V^{\sharp}(I\,\stackrel{\cdot}{\sharp}\,J)\subseteq  V^{\sharp}(I)\cup V^{\sharp}(J).
\end{equation}

\medskip
It follows from (\ref{eq50}) and (\ref{eq51}) that (ii) is true.

\hfill\raisebox{1mm}{\framebox[2mm]{}}

\bigskip
Let $(\, R=R_0\oplus R_1, \, + , \, \cdot\, , \,\sharp\, \,)$ be a Hu-Liu triring. By Proposition~\ref{pr7.1}, the collection
$$
V^{\sharp}: =\{\, V^{\sharp}(I) \, |\, \mbox{$I$ is a triideal of $R$} \, \}
$$
of subsets of $Spec^{\sharp}R$ satisfies the axioms for closed sets in a topological space. The topology on $Spec^{\sharp}R$ having the elements of $V^{\sharp}$ as closed sets is called the {\bf extended Zariski topology}. The collection
$$
D^{\sharp}: =\{\, D^{\sharp}(I) \, |\, \mbox{$I$ is a triideal of $R$} \, \}
$$
consists of the open sets of the extended Zariski topology on $Spec^{\sharp}R$, where
$$
D^{\sharp}(I): =Spec^{\sharp}R\setminus V^{\sharp}(I)=
\{\, P \, |\, \mbox{$P\in Spec^{\sharp}R$ and $P\not\supseteq I$} \, \}.
$$

\medskip
For $x_i\in R_i$ with $i\in \{0,\, 1\}$, both $Rx_0$ and $R_1\,\sharp\,x_1$ are triideals of $R$ by Proposition~\ref{pr2.1}. Let
$$
D^{\sharp}(x_0): =D^{\sharp}(Rx_0), \qquad D^{\sharp}(x_1): =D^{\sharp}(R_1\,\sharp\,x_1).
$$
If $I_0$ and $I_1$ are the even part and odd part of an triideal $I$, then 
$$D^{\sharp}(I)=\bigcup _{x_i\in I_i,\, i=0,1}D^{\sharp}(x_i).$$
Thus, $\{\, D^{\sharp}(x_i)\,|\,\mbox{$x_i\in R_i$ with  $i=0,\, 1$}\,\}$ forms an open base for the extended Zariski topology on $Spec^{\sharp}R$. Each $D^{\sharp}(x_i)$ is called a {\bf basic open subset} of $Spec^{\sharp}R$. Clearly, $D^{\sharp}(0)=\emptyset$, $D^{\sharp}(1)=Spec^{\sharp}R$, $D^{\sharp}(1^{\sharp})=Spec^{\sharp}_1R$, and $D^{\sharp}(x_1)\subseteq  Spec^{\sharp}_1R$ for $x_1\in R_1$.

\begin{proposition}\label{pr7.2} Let $I$ and $J$ be triideals of a Hu-Liu triring. 
\begin{description}
\item[(i)] $V^{\sharp}(I)\subseteq V^{\sharp}(J)$ if and only if 
$\sqrt[\sharp]{J}\subseteq \sqrt[\sharp]{I}$.
\item[(ii)] $V^{\sharp}(I)=V^{\sharp}(\sqrt[\sharp]{I})$.
\end{description}
\end{proposition}

\medskip
\noindent
{\bf Proof} (i) If $V^{\sharp}(I)\subseteq V^{\sharp}(J)$, then $P\supseteq I$ implies that $P\supseteq J$ for $P\in Spec^{\sharp}R$. Thus, we have
\begin{eqnarray*}
&&\{\,P\,|\,\mbox{$P\in Spec^{\sharp}R$ and $P\supseteq J$}\,\}\\
&=&\{\,P\,|\,\mbox{$P\in Spec^{\sharp}R$ and $P\supseteq I$}\,\}\bigcup
\{\,P\,|\,\mbox{$P\in Spec^{\sharp}R$, $P\supseteq J$ and $P\not\supseteq I$}\,\}.
\end{eqnarray*}

By Proposition~\ref{pr6.4}, we get
\begin{eqnarray*}
\sqrt[\sharp]{J}&=&\bigcap _{\begin{array}{c}P\in spec^{\sharp}R\\P\supseteq J\end{array}} P=
\Bigg(\bigcap _{\begin{array}{c}P\in spec^{\sharp}R\\P\supseteq I\end{array}} P\Bigg)\bigcap
\Bigg(\bigcap _{\begin{array}{c}P\in spec^{\sharp}R\\P\supseteq J,\,P\not\supseteq I\end{array}} P\Bigg)\\
&\subseteq &\bigcap _{\begin{array}{c}P\in spec^{\sharp}R\\P\supseteq I\end{array}} P=
\sqrt[\sharp]{J}.
\end{eqnarray*}

\medskip
Conversely, if $\sqrt[\sharp]{J}\subseteq \sqrt[\sharp]{I}$, then for any 
$Q\in V^{\sharp}(I)$, we have $Q\supseteq I$ and
$$
Q\supseteq \bigcap _{\begin{array}{c}P\in Spec^{\sharp}R\\P\supseteq I\end{array}} P
=\sqrt[\sharp]{I}\supseteq \sqrt[\sharp]{J}=
\bigcap _{\begin{array}{c}P\in Spec^{\sharp}R\\P\supseteq J\end{array}} P\supseteq J,
$$
which proves that $Q\in V^{\sharp}(J)$. Thus, we get $V^{\sharp}(I)\subseteq V^{\sharp}(J)$.

\medskip
(ii) Since $\sqrt[\sharp]{\sqrt[\sharp]{I}}=\sqrt[\sharp]{I}$, (ii) follows from (i).

\hfill\raisebox{1mm}{\framebox[2mm]{}}

\bigskip
\begin{definition}\label{def7.1} Let $X$ be a topological space.
\begin{description}
\item[(i)] A closed subset $F$ of $X$ is {\bf reducible} if $F=F_{(1)}\cup F_{(2)}$ for proper closed subsets $F_{(1)}$, $F_{(2)}$ of $X$. We call a closed subset $F$ {\bf irreducible} if it is not reducible.
\item[(ii)] $X$ is {\bf quasicompact} if given an arbitrary open covering 
$\{\, U_{(i)}\,|\, i\in I\,\}$ of $X$, there exists a finite subcovering of $X$, i.e., there exist finitely many members $U_{(i_1)}$, $\dots$, $U_{(i_n)}$ of $\{\, U_{(i)}\,|\, i\in I\,\}$ such that $X=U_{(i_1)}\cup\cdots\cup U_{(i_n)}$.
\end{description}
\end{definition}

\medskip
Using the topological concepts above, we have the following

\medskip
\begin{proposition}\label{pr7.3}  Let $(\, R=R_0\oplus R_1, \, + , \, \cdot , \, \,\sharp\, \,)$ be a Hu-Liu triring.
\begin{description}
\item[(i)] The trispectrum $Spec^{\sharp}R$ is quasicompact.
\item[(ii)] Both $Spec^{\sharp}_0R$ and $Spec^{\sharp}_1R$ are quasicompact subsets of $Spec^{\sharp}R$.
\item[(iii)] If $I$ is a triideal of $R$, then the closed subset $V^{\sharp}(I)$ of $Spec^{\sharp}R$ is irreducible if and only if $\sqrt[\sharp]{I}$ is a prime triideal.
\end{description}
\end{proposition}

\medskip
\noindent
{\bf Proof} (i) Let $\{\, D^{\sharp}(I_{(i)})\,|\, i\in\Delta\,\}$ be an open covering of $Spec^{\sharp}R$, where $I_{(i)}=I_{(i)0}\oplus I_{(i)1}$ is a triideal with the even part $I_{(i)0}$ and the odd part $I_{(i)1}$ for each $i\in \Delta$. Thus, 
$Spec^{\sharp}R=\bigcup_{i\in\Delta}D^{\sharp}(I_{(i)})=
D^{\sharp}(\sum_{i\in\Delta}I_{(i)})$. Hence, $V^{\sharp}(\sum_{i\in\Delta}I_{(i)})=\emptyset$. If $1\not\in (\sum_{i\in\Delta}I_{(i)})_0=\sum_{i\in\Delta}I_{(i)0}$, then $(\sum_{i\in\Delta}I_{(i)})_0$ is a proper ideal of the commutative ring 
$( R_0,\, +,\, \cdot)$. Hence, there exists a maximal ideal $M_0$ of the ring 
$( R_0,\, +,\, \cdot)$ such that $(\sum_{i\in\Delta}I_{(i)})_0\subseteq M_0$. Since 
$M_0\oplus R_1$ is a prime triideal of $R$ and $\sum_{i\in\Delta}I_{(i)}\subseteq M_0\oplus R_1$, we get that $M_0\oplus R_1\in V^{\sharp}(\sum_{i\in\Delta}I_{(i)})=\emptyset$, which is impossible. Therefore, $1 \in (\sum_{i\in\Delta}I_{(i)})_0=\sum_{i\in\Delta}I_{(i)0}$, which implies that $x_{(i_1)0}+x_{(i_2)0}+\cdots +x_{(i_n)0}=1$ for some positive integer $n$ and 
$x_{(i_k)0}\in I_{(i_k)0}$ with $i_k\in\Delta$ and $n\ge k\ge 1$. It follows that
$D^{\sharp}(x_{(i_k)0})\subseteq D^{\sharp}(I_{(i_k)})$ and
\begin{eqnarray*}
&&Spec^{\sharp}R=\bigcup_{i\in\Delta}D^{\sharp}(I_{(i)})\supseteq \bigcup_{k=1}^nD^{\sharp}(I_{(i_k)})\supseteq \bigcup_{k=1}^nD^{\sharp}(x_{(i_k)0})\\
&=&\bigcup_{k=1}^nD^{\sharp}(Rx_{(i_k)0})
=D^{\sharp}\left(\sum_{k=1}^nRx_{(i_k)0}\right)=D^{\sharp}(R)=Spec^{\sharp}R,
\end{eqnarray*}
which implies that $Spec^{\sharp}R=\bigcup_{k=1}^nD^{\sharp}(I_{(i_k)})$.

\bigskip
(ii) Note that a closed subset of a quasicompact topological space is a quasicompact subset. Since  $Spec^{\sharp}_0R=V^{\sharp}(R_1)$ is a closed subset of the quasicompact topological space $Spec^{\sharp}R$, $Spec^{\sharp}_0R$ is a quasicompact subset.

\medskip
It is well-known that a subset $C$ of a topological space $X$ is a quasicompact subset of $X$ if and only if every covering of $C$ by open subsets of $X$ has a finite subcovering. Hence, in order to prove that $Spec^{\sharp}_1R$ is a quasicompact subsets of $Spec^{\sharp}R$, it suffices to prove that
if $Spec^{\sharp}_1R=\bigcup_{j\in\Gamma}D^{\sharp}(J_{(j)})$ for triideals $J_{(j)}$ of $R$, then there exists a positive integer $m$ such that $Spec^{\sharp}_1R=\bigcup_{k=1}^mD^{\sharp}(J_{(j_k)})$ for some $j_1$, $\cdots$, $j_k\in\Gamma$.

\medskip
Since $Spec^{\sharp}_1R=\bigcup_{j\in\Gamma}D^{\sharp}(J_{(j)})=D^{\sharp}(\sum_{j\in\Gamma}J_{(j)})$, we have $V^{\sharp}(\sum_{j\in\Gamma}J_{(j)})=Spec^{\sharp}R\setminus D^{\sharp}(\sum_{j\in\Gamma}J_{(j)})=Spec^{\sharp}_0R$. If $(\sum_{j\in\Gamma}J_{(j)})_1\ne R_1$, then there exists a maximal ideal $N_1$ of the commutative ring 
$( R_1,\, +,\, \sharp)$ such that $(\sum_{j\in\Gamma}J_{(j)})_1\subseteq N_1$. By Proposition~\ref{pr5.1}, there exists an ideal $N_0$ of the commutative ring 
$( R_0,\, +,\, \cdot)$ such that $N_0\supseteq (\sum_{j\in\Gamma}J_{(j)})_0$ and $N_0\oplus N_1$
is a prime triideal of $R$. Thus, 
$N_0\oplus N_1\in V^{\sharp}(\sum_{j\in\Gamma}J_{(j)})=Spec^{\sharp}_0R$, which is impossible because $N_1\ne R_1$. This proves that $(\sum_{j\in\Gamma}J_{(j)})_1= R_1$. Hence, we have
$y_{(j_1)1}+y_{(j_2)1}+\cdots +y_{(j_m)1}=1^{\sharp}$ for some positive integer $m$ and $y_{(j_k)1}\in J_{(j_k)1}$ with $j_k\in\Gamma$ and $m\ge k\ge 1$. It follows that 
$D^{\sharp}(y_{(j_k)1})\subseteq D^{\sharp}(J_{(j_k)})$ and
\begin{eqnarray*}
&&Spec^{\sharp}_1R=\bigcup_{j\in\Gamma}D^{\sharp}(J_{(j)})
\supseteq \bigcup_{k=1}^mD^{\sharp}(y_{(j_k)1})
=\bigcup_{k=1}^mD^{\sharp}(R_1\,\sharp\,y_{(j_k)1})\\
&=&D^{\sharp}\left(\sum_{k=1}^m R_1\,\sharp\,y_{(j_k)1}\right)
=D^{\sharp}(R_1)=Spec^{\sharp}_1R,
\end{eqnarray*}
which implies that $Spec^{\sharp}_1R=\bigcup_{k=1}^mD^{\sharp}(y_{(j_k)1})$.

\bigskip
(iii) By Proposition~\ref{pr7.2}, we may assume $I=\sqrt[\sharp]{I}$ in the following proof.
First, we prove that if $V^{\sharp}(I)$ is irreducible , then $I$ is a prime triideal. 

\medskip
Suppose that $a_0b_0\in I$ for some $a_0$, $b_0\in R_0$. Let 
$$J_{(1)}=I+Ra_0=(I_0+R_0a_0)\oplus (I_1+R_1a_0)$$
and
$$K_{(1)}=I+Rb_0=(I_0+R_0b_0)\oplus (I_1+R_1b_0).$$
Then both $J_{(1)}$ and $K_{(1)}$ are triideals of $R$ and
\begin{eqnarray*}
J_{(1)}\,\stackrel{\cdot}{\sharp}\, K_{(1)}&=&(I_0+R_0a_0)(I_0+R_0b_0)+
(I_1+R_1a_0)\,\sharp\,(I_1+R_1b_0)\\
&\subseteq& I+R_0a_0R_0b_0+(R_1a_0)\,\sharp\,(R_1b_0)\subseteq
I+R_0a_0b_0+R_1a_0b_0\subseteq I.
\end{eqnarray*}
Hence, we get $V^{\sharp}(J_{(1)})\cup V^{\sharp}(K_{(1)})=
V^{\sharp}(J_{(1)}\,\stackrel{\cdot}{\sharp}\, K_{(1)})\supseteq V^{\sharp}(I)$ by Proposition~\ref{pr7.1} (ii). It is clear that 
$V^{\sharp}(J_{(1)})\subseteq V^{\sharp}(I)$ and $V^{\sharp}(K_{(1)})\subseteq V^{\sharp}(I)$. Hence, we get that $V^{\sharp}(J_{(1)})\cup V^{\sharp}(K_{(1)})\subseteq  V^{\sharp}(I)$. Thus we have  $V^{\sharp}(J_{(1)})\cup V^{\sharp}(K_{(1)})= V^{\sharp}(I)$. Since $V^{\sharp}(I)$ is irreducible, $V^{\sharp}(J_{(1)})=V^{\sharp}(I)$ or $V^{\sharp}(K_{(1)})=V^{\sharp}(I)$, which imply that $I= \sqrt[\sharp]{J_{(1)}}\supseteq J_{(1)}\ni a_0$ or
$I= \sqrt[\sharp]{K_{(1)}}\supseteq K_{(1)}\ni b_0$. This proves that
\begin{equation}\label{eq52} 
a_0b_0\in I\Longrightarrow \mbox{$a_0\in I$ or $b_0\in I$ for $a_0$, $b_0\in R_0$}.
\end{equation}

\medskip
Suppose that $a_0b_1\in I$ for some $a_0\in R_0$ and $b_1\in R_1$. Using the triideals 
$$J_{(2)}=I+Ra_0=(I_0+R_0a_0)\oplus (I_1+R_1a_0)$$ and 
$$K_{(2)}=I+R_1\,\sharp\, b_1=I_0\oplus (I_1+R_1\,\sharp\, b_1),$$ 
we have
\begin{eqnarray*}
J_{(2)}\,\stackrel{\cdot}{\sharp}\, K_{(2)}&=&(I_0+R_0a_0)I_0+
(I_1+R_1a_0)\,\sharp\,(I_1+R_1\,\sharp\, b_1)\\
&\subseteq& I+(R_1a_0)\,\sharp\,R_1\,\sharp\, b_1\subseteq
I+(a_0R_1)\,\sharp\,b_1\,\sharp\, R_1\subseteq I+a_0(R_1\,\sharp\,b_1)\,\sharp\, R_1\\
&\subseteq& I+a_0(b_1\,\sharp\,R_1)\,\sharp\, R_1\subseteq I+(a_0b_1)\,\sharp\,R_1\,\sharp\, R_1
\subseteq I,
\end{eqnarray*}
which implies that  $I= \sqrt[\sharp]{J_{(2)}}\supseteq J_{(2)}\ni a_0$ or
$I= \sqrt[\sharp]{K_{(2)}}\supseteq K_{(2)}\ni b_1$. This proves that
\begin{equation}\label{eq53} 
a_0b_1\in I\Longrightarrow \mbox{$a_0\in I$ or $b_1\in I$ for $a_0\in R_0$ and $b_1\in R_1$}.
\end{equation}

\medskip
Similarly, we have
\begin{equation}\label{eq54} 
a_1b_0\in I\Longrightarrow \mbox{$a_1\in I$ or $b_0\in I$ for $a_1\in R_1$ and $b_0\in R_0$}
\end{equation}
and
\begin{equation}\label{eq55} 
a_1\,\sharp\,b_1\in I\Longrightarrow \mbox{$a_1\in I$ or $b_1\in I$ for $a_1$, $b_1\in R_1$}.
\end{equation}

\medskip
By (\ref{eq52}), (\ref{eq53}), (\ref{eq54}) and (\ref{eq55}), $I$ is a prime triideal. 

\bigskip
Next, we prove that if $I$ is a prime triideal, then $V^{\sharp}(I)$ is irreducible. Suppose that 
$V^{\sharp}(I)=V^{\sharp}(J)\cup V^{\sharp}(K)$, where $J$ and $K$ are triideals of $R$. Using Proposition~\ref{pr7.2} (ii), we can assume that $\sqrt[\sharp]{J}=J$ and $\sqrt[\sharp]{K}=K$. In this case, we have
$$
V^{\sharp}(J)\subseteq V^{\sharp}(I)\Longrightarrow I=
\sqrt[\sharp]{I}\subseteq \sqrt[\sharp]{J}=J
$$
and
$$
V^{\sharp}(K)\subseteq V^{\sharp}(I)\Longrightarrow I=
\sqrt[\sharp]{I}\subseteq \sqrt[\sharp]{K}=K.
$$
By Proposition~\ref{pr7.1} (ii), we have 
$V^{\sharp}(I)=V^{\sharp}(J)\cup V^{\sharp}(K)=V^{\sharp}(J\,\stackrel{\cdot}{\sharp}\,K)$. This fact and Proposition~\ref{pr7.2} (i) give $J\,\stackrel{\cdot}{\sharp}\,K\subseteq \sqrt[\sharp]{I}=I$. Since $I$ is a prime triideal, we get $J\subseteq I$ or $K\subseteq I$ by 
Proposition~\ref{pr5.2}. Hence, $I=J$ or $I=K$. Thus, $V^{\sharp}(I)=V^{\sharp}(J)$ or 
$V^{\sharp}(I)=V^{\sharp}(K)$. This proves that $V^{\sharp}(I)$ is irreducible.

\hfill\raisebox{1mm}{\framebox[2mm]{}}

\bigskip
\section{Localization of Hu-Liu Trirings}

Let $(\, R=R_0\oplus R_1, \, + , \, \cdot , \, \,\sharp\, \,)$ be a Hu-Liu triring with the identity $1$ and the local identity $1^{\sharp}$. A subset $S$ of $R$ is called a 
{\bf multiplicative subset} if 
$$
0\not\in S, \quad 1\in S_0:=S\cap R_0, \quad 1^{\sharp}\in S_1:=S\cap R_1, \quad 
S=S_0\cup S_1
$$ 
and
\begin{equation}\label{eq11}
s_0t_0\in S_0, \quad s_0s_1\in S_1, \quad s_1s_0\in S_1,\quad s_1 \,\sharp\, t_1\in S_1,
\quad s_0S_1=S_1s_0,
\end{equation}
where $s_i$, $t_i\in S_i$ and $i=0$, $1$.

\medskip
Given a  multiplicative subset $S$ of a Hu-Liu triring $R$, we define a relation in the Cartesian product $(R_0\times S_0)\times (R_1\times S_1)$ as follows:
\begin{eqnarray}
&&\left((a_0, s_0),\, (a_1, s_1) \right)\sim \left((b_0, t_0),\, (b_1, t_1)\right)\nonumber\\
\label{eq12}&\Leftrightarrow& u_0(a_0t_0-b_0s_0)=0 \quad\mbox{and}\quad u_1\,\sharp\,(a_1\,\sharp\,t_1\, -b_1\,\sharp\,s_1)=0
\end{eqnarray}
for some $u_i\in S_i$, where $a_i$, $b_i\in R_i$, $s_i$, $t_i\in S_i$ and $i=0$, $1$.

It is clear that $\sim$ is an equivalence relation. Let 
$\left(\displaystyle\frac{a_0}{s_0},  \,\displaystyle\frac{a_1}{s_1}\right)$ be the equivalence class containing $\left((a_0, s_0),\, (a_1, s_1) \right)$; that is,
$$
\left(\displaystyle\frac{a_0}{s_0},  \,\displaystyle\frac{a_1}{s_1}\right):=
\left\{\, \left((a_0', s_0'),\, (a_1', s_1') \right) \, \left| \, \begin{array}{c} 
\mbox{$(a_i', s_i')\in R_i\times S_i$,}\\
\mbox{$u_0(a_0s_0'-s_0a_0')=0$,}\\
\mbox{ $u_1\,\sharp\,(a_1\,\sharp\,s_1'\, -s_1\,\sharp\,a_1')=0$ }\\
\mbox{for some $u_i\in S_i$ with $i=0$, $1$.}
\end{array}\right.\right\}.
$$

\medskip
Let
$$
S^{\stackrel{\sharp}{-}1}R:=\left\{\, \left.\left(\displaystyle\frac{a_0}{s_0},  \,\displaystyle\frac{a_1}{s_1}\right) \, \right | \, 
\left((a_0, s_0),\, (a_1, s_1)\right)\in (R_0\times S_0)\times (R_1\times S_1)\right\}
$$
be the set of all equivalence classes. Then 
$S^{\stackrel{\sharp}{-}1}R=(S^{\stackrel{\sharp}{-}1}R)_0\oplus 
(S^{\stackrel{\sharp}{-}1}R)_1$ is a Hu-Liu triring, where the even part $\Big(S^{\stackrel{\sharp}{-}1}R\Big)_0$ and odd part 
$\Big(S^{\stackrel{\sharp}{-}1}R\Big)_1$ are given by
$$
\Big(S^{\stackrel{\sharp}{-}1}R\Big)_0:=\left\{\left.\left(\displaystyle\frac{a_0}{s_0}, \,\displaystyle\frac{0}{1^\sharp}\right) \, \right| \, (a_0, s_0)\in R_0\times S_0 \,\right\},
$$
$$
\Big(S^{\stackrel{\sharp}{-}1}R\Big)_1:=\left\{\left.\left(\displaystyle\frac{0}{1}, \,\displaystyle\frac{a_1}{s_1}\right) \, \right| \, (a_1, s_1)\in R_1\times S_1 \,\right\}
$$
and the addition $+$, the ring multiplication $\cdot$ and the local product $\,\sharp\,$ are given by
$$
\left(\displaystyle\frac{a_0}{s_0},  \,\displaystyle\frac{a_1}{s_1}\right)+
\left(\displaystyle\frac{b_0}{t_0},  \,\displaystyle\frac{b_1}{t_1}\right):=
\left(\displaystyle\frac{a_0t_0+s_0b_0}{s_0t_0},  
\,\displaystyle\frac{a_1\,\sharp\, t_1+s_1\,\sharp\, b_1}{s_1\,\sharp\, t_1}\right),
$$
$$
\left(\displaystyle\frac{a_0}{s_0},  \,\displaystyle\frac{a_1}{s_1}\right)\cdot
\left(\displaystyle\frac{b_0}{t_0},  \,\displaystyle\frac{b_1}{t_1}\right):=
\left(\displaystyle\frac{a_0b_0}{s_0t_0},  \,\displaystyle\frac{(a_0b_1)\,\sharp\,(s_1t_0)+(s_0t_1)\,\sharp\,(a_1b_0)}{(s_0t_1)\,\sharp\, (s_1t_0)}\right)
$$
and
$$
\left(\displaystyle\frac{0}{1}, \,\displaystyle\frac{a_1}{s_1}\right) \,\sharp\,
\left(\displaystyle\frac{0}{1}, \,\displaystyle\frac{b_1}{t_1}\right):=
\left(\displaystyle\frac{0}{1}, \,\displaystyle\frac{a_1\,\sharp\,b_1}{s_1\,\sharp\,t_1}\right).
$$
Clearly, $\left(\displaystyle\frac{1}{1}, \,\displaystyle\frac{0}{1^{\sharp}}\right)$ is the identity of 
the ring $S^{\stackrel{\sharp}{-}1}R$, and 
$\left(\displaystyle\frac{0}{1},  \,\displaystyle\frac{1^{\sharp}}{1^{\sharp}}\right)$ is the local identity.

\medskip
Let $i^{\sharp S}_R: R\to S^{\stackrel{\sharp}{-}1}R$ be the map defined by
$$
i^{\sharp S}_R: a=a_0+a_1\mapsto \left(\displaystyle\frac{a_0}{1}, \,\displaystyle\frac{a_1}{1^{\sharp}}\right),
$$
where $a_i\in R_i$ for $i=0$ and $1$. Then $i^{\sharp S}_R$ is a triring homomorphism from $R$ to 
$S^{\stackrel{\sharp}{-}1}R$.

\medskip
The pair $\Big(S^{\stackrel{\sharp}{-}1}R, \, i^{\sharp S}_R \Big)$ constructed above is called the {\bf localization} of $R$ with respect to $S$, and $i^{\sharp S}_R$ is called the 
{\bf canonical triring homomorphism} from $R$ to $S^{\stackrel{\sharp}{-}1}R$. The next proposition gives the universal mapping characterization of the localization 
$\Big(S^{\stackrel{\sharp}{-}1}R, \, i^{\sharp S}_R \Big)$.

\begin{proposition}\label{pr3.1} Let $R=R_0\oplus R_1$ be a Hu-Liu triring. If 
$S=S_0\cup S_1$ is a multiplicative subset of $R$, then the 
localization $\Big(S^{\stackrel{\sharp}{-}1}R, \, i^{\sharp S}_R \Big)$ of $R$ with respect to $S$ has the following two properties.
\begin{description}
\item[(i)] For any $s_i\in S_i$ with $i=0$ and $1$, $i^{\sharp S}_R(s_0)$ is invertible in the commutative ring $\left(\Big(S^{\stackrel{\sharp}{-}1}R\Big)_0, +, \cdot \right)$ and $i^{\sharp S}_R(s_1)$ is invertible in the commutative ring $\left(\Big(S^{\stackrel{\sharp}{-}1}R\Big)_1, +, \,\sharp\, \right)$.
\item[(ii)] If $\psi: R\to \overline{R}$ is a triring homomorphism from the Hu-Liu triring $R$ to a Hu-Liu triring $\overline{R}$ such that $\psi (s_0)$ is invertible in the commutative ring 
$(\overline{R}_0, +, \cdot )$ and $\psi (s_1)$ is invertible in the commutative ring 
$(\overline{R}_1, +, \,\sharp\, )$ for any $s_i\in S_i$ with $i=0$ and $1$, then there exists a unique triring homomorphism
$\overline{\psi}: S^{\stackrel{\sharp}{-}1}R\to \overline{R}$ such that $\psi=\overline{\psi}\,i^{\sharp S}_R$.
\end{description}
\end{proposition}

\medskip
\noindent
{\bf Proof} A direct computation.

\hfill\raisebox{1mm}{\framebox[2mm]{}}

\bigskip
Let  $R=R_0\oplus R_1$ be a Hu-Liu triring. The following three types of localizations of $R$ are very useful in the study of Hu-Liu trirings: 
\begin{description}
\item[Type 1.] If $P=P_0\oplus P_1$ is a triideal of $R$ with $P_1\ne R_1$, then $(R_0\setminus P_0)\bigcup (R_1\setminus P_1)$ is a multiplicative subset of $R$. The localization of $R$ with respect to $(R_0\setminus P_0)\cup (R_1\setminus P_1)$ is called the {\bf localization} of $R$ at $P$ and is denoted by $R^{\sharp}_P$. 
\item[Type 2.] If $f_0\in R_0$ and $D^{\sharp}(f_0)\bigcap spec^{\sharp}_1R\ne\emptyset$, then
$T(f_0)=T_0(f_0)\cup T_1(f_0)$ is a multiplicative subset of $R$, where 
$T_0(f_0):=\{\, f_0^n\,|\, n\in\mathcal{Z}_{\ge 0}\,\}$
and
$$T_1(f_0):=\{\, z_1f_0^n\,|\,\mbox{ $n\in\mathcal{Z}_{\ge 0}$ and $z_1\in R_1\setminus P$ for all $P\in D^{\sharp}(f_0)\bigcap spec^{\sharp}_1R$}\,\}.$$
The localization of $R$ with respect to $T(f_0)=T_0(f_0)\cup T_1(f_0)$ is called the {\bf localization of $R$ at $f_0$} and is denoted by $R_{f_0}^{\sharp}$. Clearly, we have
$$
R_{f_0}^{\sharp}=\left\{\,\left.\frac{a_0}{f_0^n}+\frac{a_1}{z_1f_0^m}\,\right|\, \mbox{$m$, $n\ge 0$, $a_i\in R_i$ for $i=0$, $1$ and $z_1f_0^m\in T_1(f_0)$}\,\right\}.
$$
\item[Type 3.] If $f_1\in R_1$ and $f_1$ is not trinilpotent, then 
$\{\, 1\}\bigcup \{\, f_1^{\sharp\,m}\,\}_{m\ge 0}$
is multiplication subset of $R$. The localization of $R$ with respect to 
$\{\, 1\}\bigcup \{\, f_1^{\sharp\,m}\,\}_{m\ge 0}$ is called the {\bf localization of $R$ at $f_1$} and is denoted by $R_{f_1}^{\sharp}$. Thus, we have

$$
R_{f_1}^{\sharp}=\left\{\,\left.a_0+\frac{a_1}{f_1^{\sharp\,m}}\,\right|\, \mbox{$m$, $n\ge 0$ and $a_i\in R_i$ for $i=0$, $1$}\,\right\}.
$$
\end{description}

\bigskip
\section{Sheaf Structures on Trispectra}

We begin this section with the following

\begin{definition}\label{def8.1} Let $X$ be a topological space. A {\bf presheaf} $\mathcal{P}$ of Hu-Liu trirings on $X$ assigns a Hu-Liu triring $\mathcal{P}(U)$ to each open set $U$ of $X$ and a triring homomorphism: $\rho _{_{UV}}: \mathcal{P}(U)\to \mathcal{P}(V)$, called the {\bf restriction map}, to each inclusion of open sets $V\subseteq U$, subject to the conditions:
\begin{description}
\item[(i)] $\mathcal{P}(0)=\emptyset$;
\item[(ii)] $\rho _{_{UU}}: \mathcal{P}(U)\to \mathcal{P}(U)$ is the identity map;
\item[(iii)] If $W\subseteq V\subseteq U$ are three open subsets, then 
$\rho _{_{UW}}=\rho _{_{VW}}\rho _{_{UV}}$.
\end{description}
\end{definition}

\medskip
If $\mathcal{P}$ is a presheaf of Hu-Liu trirings on a  topological space $X$, then an element $s\in \mathcal{P}(U)$ is called a {\bf section} of $\mathcal{P}$ over the open set $U$.

\medskip
\begin{definition}\label{def8.2} A presheaf $\mathcal{P}$ of Hu-Liu trirings on a  topological space $X$ is called a {\bf sheaf} of  Hu-Liu trirings if it satisfies the following two conditions for an arbitrary open subset $U$ of $X$ and an arbitrary open covering 
$\{\,U_{(i)}\,|\, i\in\Delta\,\}$ of $U$:
\begin{description}
\item[(i)] If $s\in \mathcal{P}(U)$ and $\rho _{_{UU_{(i)}}}(s)=0$ for all $i\in\Delta$, then $s=0$;
\item[(ii)] If we have elements $s_{(i)}\in \mathcal{P}(U_{(i)})$ for each $i\in\Delta$ having the property that 
$\rho _{_{U_{(i)}, U_{(i)}\cap U_{(j)}}}(s_{(i)})=
\rho _{_{U_{(j)}, U_{(i)}\cap U_{(j)}}}(s_{(j)})$ for $i$, $j\in \Delta$, then there is an element $s\in \mathcal{P}(U)$ such that $\rho _{_{UU_{(i)}}}(s)=s_{(i)}$ for each $i\in\Delta$.
\end{description}
\end{definition}

\medskip
Let $(\, R=R_0\oplus R_1, \, + , \, \cdot , \, \,\sharp\, \,)$ be a Hu-Liu triring with the identity $1$ and the local identity $1^{\sharp}$. If $f_0\in R_0$ and $f_1\in R_1$, then we define
$$
\mathcal{O}(D^{\sharp}(f_0)):=\left\{\begin{array}{ll}
(R_0)_{f_0}\oplus R_0\, &\,\mbox{if $f_0\not\in \sqrt[\sharp]{0}$ and $D^{\sharp}(f_0)\bigcap spec^{\sharp}_1R=\emptyset$},\\
R_{f_0}^{\sharp}\oplus R_{1^{\sharp}f_0}^{\sharp}\, &\,\mbox{if $f_0\not\in \sqrt[\sharp]{0}$ and $D^{\sharp}(f_0)\bigcap spec^{\sharp}_1R\ne\emptyset$},\\
0\, &\,\mbox{if $f_0\in \sqrt[\sharp]{0}$}\end{array}\right.
$$
and
$$
\mathcal{O}(D^{\sharp}(f_1)):=\left\{\begin{array}{ll}R_{f_1}^{\sharp}\, &\,\mbox{if $f_1\not\in \sqrt[\sharp]{0}$},\\
0\, &\,\mbox{if $f_1\in \sqrt[\sharp]{0}$}\end{array}\right. ,
$$
where $(R_0)_{f_0}$ is the ordinary localization of the commutative ring $(R_0, \, +,\, \cdot)$ with respect to the multiplication set $\{\,f^n_0\,|\, n\in\mathcal{Z}_{\ge 0}\,\}$, and 
$(R_0)_{f_0}\oplus R_0$ is a Hu-Liu triring with zero odd part. We then have the following result.

\medskip
\begin{proposition}\label{pr8.1} Let $(\, R=R_0\oplus R_1, \, + , \, \cdot , \, \,\sharp\, \,)$ be a Hu-Liu triring. If $f_i$, $g_i\in R_i$  for $i=0$, $1$ and 
$D^{\sharp}(g_i)\subseteq D^{\sharp}(f_j)$ for $i$, $j\in\{\, 0,\, 1\,\}$, then there exist a triring homomorphism $\rho _{_{f_j,g_i}}: \mathcal{O}(D^{\sharp}(f_j))\to \mathcal{O}(D^{\sharp}(g_i))$ such that
\begin{description}
\item[(i)] $\rho _{_{f_j, f_j}}=id_{_{D^{\sharp}(f_j)}}$,
\item[(ii)] $\rho _{_{f_j, h_k}}=\rho _{_{g_i, h_k}}\rho _{_{f_j, g_i}}$ provided 
$D^{\sharp}(h_k)\subseteq D^{\sharp}(g_i)\subseteq D^{\sharp}(f_j)$,
\end{description}
where $i$, $j$, $k\in\{\, 0, \, 1\,\}$.
\end{proposition}

\medskip
\noindent
{\bf Proof} If $g_i\in \sqrt[\sharp]{0}$, then we define $\rho _{_{f_j,g_i}}: \mathcal{O}(D^{\sharp}(f_j))\to \mathcal{O}(D^{\sharp}(g_i))$ to be the zero homomorphism, i.e.,
\begin{equation}\label{eq56} 
\rho _{_{f_j,g_i}}(x):=0\quad\mbox{for $x\in \mathcal{O}(D^{\sharp}(f_j))$ and 
$g_i\in \sqrt[\sharp]{0}$}.
\end{equation}

\medskip
We now define $\rho _{_{f_j,g_i}}: \mathcal{O}(D^{\sharp}(f_j))\to \mathcal{O}(D^{\sharp}(g_i))$ for $g_i\not\in \sqrt[\sharp]{0}$ and $f_j\not\in \sqrt[\sharp]{0}$ by cases.

\medskip
\underline{\it Case 1}: $i=0$, $j=0$, in which case, we have
\begin{eqnarray}\label{eq57}
&&D^{\sharp}(Rg_0)=D^{\sharp}(g_0)\subseteq D^{\sharp}(f_0)=D^{\sharp}(Rf_0)\nonumber\\
&\Longrightarrow& V^{\sharp}(Rg_0)\supseteq V^{\sharp}(Rf_0)
\Longrightarrow \sqrt[\sharp]{Rg_0}\subseteq \sqrt[\sharp]{Rf_0}
\quad\mbox{(by Proposition~\ref{pr7.2} (i))}\nonumber\\
&\Longrightarrow& g_0^u=r_0f_0\quad\mbox{for some $u\in\mathcal{Z}_{>0}$ and $r_0\in R_0$.}
\end{eqnarray}

Since $D^{\sharp}(g_0)\bigcap spec^{\sharp}_1R\subseteq 
D^{\sharp}(f_0)\bigcap spec^{\sharp}_1R$, we have three subcases. 

\medskip
\underline{\it Case 1(i)}: $D^{\sharp}(g_0)\bigcap spec^{\sharp}_1R=\emptyset$ and
$D^{\sharp}(f_0)\bigcap spec^{\sharp}_1R=\emptyset$. In this case,
using (\ref{eq57}) and universal property of the ordinary localization $(R_0)_{f_0}$, we  get a triring homomorphism 
$$\rho _{_{f_0,g_0}}: \mathcal{O}(D^{\sharp}(f_0))=(R_0)_{f_0}\oplus R_0\to (R_0)_{g_0}\oplus R_0=\mathcal{O}(D^{\sharp}(g_0))$$
such that
\begin{equation}\label{eq591} 
\rho _{_{f_0,g_0}}\left(\frac{a_0}{f_0^n},\, b_0\right)
=\left(\frac{a_0r_0^n}{g_0^{nu}},\, b_0\right)
\qquad\mbox{for $a_0,\, b_0\in R_0$}.
\end{equation}

\medskip
\underline{\it Case 1(ii)}: $D^{\sharp}(g_0)\bigcap spec^{\sharp}_1R=\emptyset$ and
$D^{\sharp}(f_0)\bigcap spec^{\sharp}_1R\ne\emptyset$. In this case, we have a triring homomorphism 
$$\rho _{_{f_0,g_0}}: \mathcal{O}(D^{\sharp}(f_0))=R_{f_0}^{\sharp}\oplus R_{1^{\sharp}f_0}^{\sharp}\to (R_0)_{g_0}\oplus R_0=\mathcal{O}(D^{\sharp}(g_0))$$
such that
\begin{equation}\label{eq592} 
\rho _{_{f_0,g_0}}\left(\left(\frac{a_0}{f_0^n}+\frac{a_1}{z_1f_0^m},\,\,
b_0+\frac{b_1}{1^{\sharp}f_0^k}\right)\right)
=\left(\frac{a_0r_0^n}{g_0^{nu}},\, b_0\right),
\end{equation}
where $z_1f_0^m\in T_1(f_0)$, $n$, $m$, $k\in \mathcal{Z}_{\ge 0}$, $a_i$, $b_i\in R_i$ and $i=0$, $1$.

\medskip
\underline{\it Case 1(iii)}: $D^{\sharp}(g_0)\bigcap spec^{\sharp}_1R\ne\emptyset$ and
$D^{\sharp}(f_0)\bigcap spec^{\sharp}_1R\ne\emptyset$.
Using (\ref{eq57}) and universal property of the localization $R_{f_0}^{\sharp}$ and $R_{1^{\sharp}f_0}^{\sharp}$, 
we get two triring homomorphism $\overline{\phi}: R_{f_0}^{\sharp}\to R_{g_0}^{\sharp}$ and
$\overline{\psi}: R_{1^{\sharp}f_0}^{\sharp}\to R_{1^{\sharp}g_0}^{\sharp}$
such that
\begin{equation}\label{eq58} 
\overline{\phi}\left(\frac{a_0}{f_0^n}+\frac{a_1}{z_1f_0^m}\right)=
\frac{a_0r_0^n}{g_0^{nu}}+\frac{a_1r_0^m}{z_1g_0^{mu}}
\end{equation}
and
\begin{equation}\label{eq59} 
\overline{\psi}\left(b_0+\frac{b_1}{1^{\sharp}f_0^k}\right)=
b_0+\frac{b_1r_0^k}{1^{\sharp}g_0^{ku}},
\end{equation}
where $z_1f_0^m\in T_1(f_0)$, $n$, $m$, $k\in \mathcal{Z}_{\ge 0}$, $a_i$, $b_i\in R_i$ and $i=0$, $1$.

By (\ref{eq58}) and (\ref{eq59}), we get a triring homomorphism 
$$\rho _{_{f_0,g_0}}: \mathcal{O}(D^{\sharp}(f_0))=R_{f_0}^{\sharp}\oplus R_{1^{\sharp}f_0}^{\sharp}\to 
R_{g_0}^{\sharp}\oplus R_{1^{\sharp}g_0}^{\sharp}=\mathcal{O}(D^{\sharp}(g_0))$$ 
such that
\begin{eqnarray}\label{eq60}
&&\rho _{_{f_0,g_0}}\left(\left(\frac{a_0}{f_0^n}+\frac{a_1}{z_1f_0^m},\,\,
b_0+\frac{b_1}{1^{\sharp}f_0^k}\right)\right)\nonumber\\
&=&\left(\overline{\phi}\left(\frac{a_0}{f_0^n}+\frac{a_1}{z_1f_0^m}\right),\,\,
\overline{\psi}\left(b_0+\frac{b_1}{1^{\sharp}f_0^k}\right)\right)\nonumber\\
&=&\left(\frac{a_0r_0^n}{g_0^{nu}}+\frac{a_1r_0^m}{z_1g_0^{mu}},\,\,
b_0+\frac{b_1r_0^k}{1^{\sharp}g_0^{ku}}\right),
\end{eqnarray}
where $z_1f_0^m\in T_1(f_0)$, $n$, $m$, $k\in \mathcal{Z}_{\ge 0}$, $a_i$, $b_i\in R_i$ and $i=0$, $1$.

\medskip
\underline{\it Case 2}: $i=0$ and $j=1$, in which case, we have 
\begin{equation}\label{eq61} 
D^{\sharp}(g_0)\subseteq D^{\sharp}(f_1)\Longrightarrow g_0\in \sqrt[\sharp]{0}. 
\end{equation}
Hence, $\rho _{_{f_1,g_0}}=0$ in this case according to the definition given at the beginning.

\medskip
\underline{\it Case 3}: $i=1$ and $j=0$, in which case, 
\begin{eqnarray}\label{eq62}
&&D^{\sharp}(g_1)\subseteq D^{\sharp}(f_0)\nonumber\\
&\Longrightarrow& g_1^{\sharp v}=r_1f_0
\quad\mbox{for some $v\in\mathcal{Z}_{>0}$ and $r_1\in R_1$}
\end{eqnarray}
and $D^{\sharp}(f_0)\bigcap spec^{\sharp}_1R\supseteq D^{\sharp}(g_1)\ne\emptyset$
It follows that we have a triring homomorphism 
$$\rho _{_{f_0,g_1}}: \mathcal{O}(D^{\sharp}(f_0))
=R_{f_0}^{\sharp}\oplus R_{1^{\sharp}f_0}^{\sharp}\to 
R_{g_1}^{\sharp}=\mathcal{O}(D^{\sharp}(g_1))$$ 
such that
\begin{equation}\label{eq63}
\rho _{_{f_0,g_0}}\left(\left(\frac{a_0}{f_0^n}+\frac{a_1}{z_1f_0^m},\,\,
b_0+\frac{b_1}{1^{\sharp}f_0^k}\right)\right)
=b_0+\frac{b_1\,\sharp\,r_1^{\sharp k}}{g_1^{\sharp kv}},
\end{equation}
where $z_1f_0^m\in T_1(f_0)$, $n$, $m$, $k\in \mathcal{Z}_{\ge 0}$, $a_i$, $b_i\in R_i$ and $i=0$, $1$.

\medskip
\underline{\it Case 4}: $i=1$ and $j=1$, in which case, 
\begin{eqnarray}\label{eq64}
&&D^{\sharp}(g_1)\subseteq D^{\sharp}(f_1)\nonumber\\
&\Longrightarrow& g_1^{\sharp \theta}=t_1\,\sharp\, f_1
\quad\mbox{for some $\theta\in\mathcal{Z}_{>0}$ and $t_1\in R_1$.}
\end{eqnarray}
It follows that we have a triring homomorphism 
$$\rho _{_{f_1,g_1}}: \mathcal{O}(D^{\sharp}(f_1))=R_{f_1}^{\sharp}
\to R_{g_1}^{\sharp}=\mathcal{O}(D^{\sharp}(g_1))$$ 
such that
\begin{equation}\label{eq65}
\rho _{_{f_1,g_1}}\left(b_0+\frac{b_1}{f_1^{\sharp n}}\right)
=b_0+\frac{b_1\,\sharp\,t_1^{\sharp n}}{g_1^{\sharp \theta n}},
\end{equation}
where $n\in \mathcal{Z}_{\ge 0}$, $b_i\in R_i$ and $i=0$, $1$.

\medskip
It is clear that the triring homomorphisms $\rho _{_{f_j,g_i}}$ defined by (\ref{eq56}), (\ref{eq60}), (\ref{eq63}) and (\ref{eq65}) has the property (i).

\medskip
Note that if $h_k\in \sqrt[\sharp]{0}$, then 
the property (ii) holds because $\rho _{_{f_j, h_k}}=\rho _{_{g_i, h_k}}=0$. 
If $h_k\not\in \sqrt[\sharp]{0}$, then we have $k\ge i\ge j$ by (\ref{eq61}). Therefore, there are only four cases:
$$
(k, i, j)=(0, 0, 0)\,\, \mbox{or}\,\, (1, 0, 0)\,\, \mbox{or}\,\, (1, 1, 0)\,\, 
\mbox{or}\, (1, 1, 1).
$$
One can check that  the property (ii) also holds  for each of the four cases above.

\hfill\raisebox{1mm}{\framebox[2mm]{}}

\bigskip
Let $U$ be a nonempty open subset of $Spec^{\sharp}R$, and let
$$S_U:=\{\, D^{\sharp}(f_{(\alpha)}) \,|\, \emptyset\ne D^{\sharp}(f_{(\alpha)})\subseteq U\quad\mbox{and}\quad \alpha \in \Lambda _U\,\}$$
be the set of the nonempty basic open subsets contained in $U$. The set $\Lambda _U$ is a partial order set with respect to the following partial order:
$$
\alpha \ge \beta\quad\mbox{if and only if}\quad D^{\sharp}(f_{(\alpha)})\supseteq D^{\sharp}(f_{(\beta)}).
$$
$\left(\mathcal{O}(D^{\sharp}(f_{(\alpha)})), \, \rho _{_{f_{(\alpha)}, f_{(\beta)}}}\right)$ is an 
{\bf inverse system} on the set $\Lambda _U$, where the triring homomorphism 
$$\rho _{_{f_{(\alpha)}, f_{(\beta)}}}: \mathcal{O}(D^{\sharp}(f_{(\alpha)}))\to \mathcal{O}(D^{\sharp}(f_{(\beta)}))\quad\mbox{for $\alpha \ge \beta$}
$$
is defined by Proposition~\ref{pr8.1}. We define
\begin{equation}\label{eq66}
\mathcal{O}(U):=\lim_{\longleftarrow _{\Lambda _U}} \mathcal{O}(D^{\sharp}(f_{(\alpha)})),
\end{equation}
where $\displaystyle\lim_{\longleftarrow _{\Lambda _U}} \mathcal{O}(D^{\sharp}(f_{(\alpha)}))$ is the {\bf inverse limit} of the inverse system \newline
$\left(\mathcal{O}(D^{\sharp}(f_{(\alpha)})), \, \rho _{_{f_{(\alpha)}, f_{(\beta)}}}\right)$.
The inverse limit is a subtriring of the direct product 
$\displaystyle\prod _{\alpha \in \Lambda _U}\mathcal{O}(D^{\sharp}(f_{(\alpha)}))$, which is given by
\begin{equation}\label{eq67}
\displaystyle\lim_{\longleftarrow _{\Lambda _U}} \mathcal{O}(D^{\sharp}(f_{(\alpha)}))=
\left\{\,(x_{(\alpha)})_{\alpha \in \Lambda _U}\,\left| \begin{array}{c}
x_{(\alpha)}\in\mathcal{O}(D^{\sharp}(f_{(\alpha)}))\\\mbox{and 
$x_{(\beta)}=\rho _{_{f_{(\alpha)}, f_{(\beta)}}}(x_{(\alpha)})$}\\ \mbox{whenever 
$\alpha \ge \beta$}\end{array}
\right.\right\}.
\end{equation}
The even part and odd part of 
$\displaystyle\lim_{\longleftarrow _{\Lambda _U}} \mathcal{O}(D^{\sharp}(f_{(\alpha)}))$ are given by
\begin{equation}\label{eq68}
\left(\displaystyle\lim_{\longleftarrow _{\Lambda _U}} \mathcal{O}(D^{\sharp}(f_{(\alpha)}))\right)_i=
\left\{\,(x_{(\alpha)i})_{\alpha \in \Lambda _U}\,\left| \begin{array}{c}
x_{(\alpha)i}\in\left(\mathcal{O}(D^{\sharp}(f_{(\alpha)}))\right)_i\\\mbox{and 
$x_{(\beta)i}=\rho _{_{f_{(\alpha)}, f_{(\beta)}}}(x_{(\alpha)i})$}\\ \mbox{whenever 
$\alpha \ge \beta$}\end{array}
\right.\right\},
\end{equation}
where $i=0$, $1$. For each $\alpha \in \Lambda _U$, let 
$p^U_{(\alpha)}: \displaystyle\lim_{\longleftarrow _{\Lambda _U}} \mathcal{O}(D^{\sharp}(f_{(\alpha)}))\to \mathcal{O}(D^{\sharp}(f_{(\alpha)}))$ be the triring homomorphism defined by
\begin{equation}\label{eq69}
p^U_{(\alpha)}\left((x_{(\gamma)})_{\gamma \in \Lambda _U}\right):=x_{(\alpha)}.
\end{equation}
Clearly, we have
\begin{equation}\label{eq70}
p^U_{(\beta)}=\rho _{_{f_{(\alpha)}, f_{(\beta)}}}p^U_{(\alpha)}\quad\mbox{if $\alpha \ge \beta$ and $\alpha$, $\beta\in \Lambda _U$ }.
\end{equation}
The pair $\left(\displaystyle\lim_{\longleftarrow _{\Lambda _U}} \mathcal{O}(D^{\sharp}(f_{(\alpha)})),\, \{\,p^U_{(\alpha)}\,|\, \alpha\in \Lambda _U\,\}\right)$
has the following universal property:
\begin{quote}
Let $(\,X,\, \{\,q_{(\alpha)}\,|\, \alpha\in \Lambda _U\,\})$ be a pair consisting of a triring $X$ and a family of triring homomorphism 
$$\{\,q_{(\alpha)}: X\to \mathcal{O}(D^{\sharp}(f_{(\alpha)}))\,|\,\alpha\in \Lambda _U \,\}.$$ If 
$q_{(\beta)}=\rho _{_{f_{(\alpha)}, f_{(\beta)}}}q_{(\alpha)}$ whenever $\alpha \ge \beta$, then there exists a unique triring homomorphism $\sigma : X\to  \displaystyle\lim_{\longleftarrow _{\Lambda _U}} \mathcal{O}(D^{\sharp}(f_{(\alpha)}))$ such that $q_{(\alpha)}=p^U_{(\alpha)}\sigma$ for each $ \alpha\in \Lambda _U$.
\end{quote}
If $V\subseteq U$, then $S_V\subseteq S_U$ and $\Lambda _V\subseteq \Lambda _U$. Consider the pair 
$$
\left(\displaystyle\lim_{\longleftarrow _{\Lambda _U}} \mathcal{O}(D^{\sharp}(f_{(\alpha)})),\, \{\,p^U_{(\alpha)}\,|\, 
\alpha\in \Lambda _V\,\}\right).
$$
Note that (\ref{eq70}) is clearly true for $\alpha$, 
$\beta\in \Lambda _V\subseteq \Lambda _U$ with $\alpha \ge \beta$. By the universal property of the following pair
$$
\left(\displaystyle\lim_{\longleftarrow _{\Lambda _V}} \mathcal{O}(D^{\sharp}(f_{(\alpha)})),\, \{\,p^U_{(\alpha)}\,|\, 
\alpha\in \Lambda _V\,\}\right),
$$
we have a unique triring homomorphism
\begin{equation}\label{eq71}
\rho _{_{U, V}}: \mathcal{O}(U)=\lim_{\longleftarrow _{\Lambda _U}} \mathcal{O}(D^{\sharp}(f_{(\alpha)}))\to
\lim_{\longleftarrow _{\Lambda _V}} \mathcal{O}(D^{\sharp}(f_{(\alpha)}))=\mathcal{O}(V)
\end{equation}
such that 
\begin{equation}\label{eq72}
p^U_{(\alpha)}=p^V_{(\alpha)}\rho _{_{U, V}}\qquad\mbox{for $\alpha\in \Lambda _V$.}
\end{equation}
It follows from (\ref{eq71}) and (\ref{eq72}) that
\begin{equation}\label{eq73}
\rho _{_{U, U}}=id_{\mathcal{O}(U)}\qquad\mbox{for each open subset $U$ of $Spec^{\sharp}_R$}
\end{equation}
and
\begin{equation}\label{eq74}
\rho _{_{U, W}}=\rho _{_{V, W}}\rho _{_{U, V}}\qquad\mbox{for three open subsets 
$W\subseteq V\subseteq U$.}
\end{equation}

By (\ref{eq73}) and (\ref{eq74}), $\mathcal{O}$ is a presheaf with the restriction map 
$\rho _{_{U, V}}$ given by (\ref{eq71}). This  presheaf is called the {\bf structure presheaf}
on $Spec^{\sharp}R$.

\medskip
\begin{proposition}\label{pr8.2} If $(\, R=R_0\oplus R_1, \, + , \, \cdot , \, \,\sharp\, \,)$ is a Hu-Liu triring, then the structure presheaf on $Spec^{\sharp}R$ is a sheaf of Hu-Liu trirings.
\end{proposition}

\medskip
\noindent
{\bf Proof} According to the definition of the structure presheaf on $Spec^{\sharp}R$, it suffices to prove that the properties (i) and (ii) in Definition~\ref{def8.2} holds for an arbitrary basic open set and an arbitrary open covering which consists of basic open sets.
Let $D^{\sharp}(f_i)$ be an arbitrary basic open set of $Spec^{\sharp}R$. Consider the following 
arbitrary open covering  of $D^{\sharp}(f_i)$
\begin{equation}\label{eq75}
D^{\sharp}(f_i)=\bigg(\bigcup _{\alpha\in\Lambda _0}D^{\sharp}(g_{(\alpha)0})\bigg)\bigcup
\bigg(\bigcup _{\beta\in\Lambda _1}D^{\sharp}(g_{(\beta)1})\bigg),
\end{equation}
where $f_0$, $g_{(\alpha)0}\in R_0$ for $\alpha\in\Lambda _0$, 
$f_1$, $g_{(\beta)1}\in R_1$ for $\beta\in\Lambda _1$, $\Lambda _0$ and $\Lambda _1$ are two index sets.

\bigskip
To prove the properties (i) in Definition~\ref{def8.2}, we need to prove that if 
$s\in \mathcal{O}(D^{\sharp}(f_i))$ satisfies
\begin{equation}\label{eq76}
\rho _{_{f_{i}, g_{(\alpha)0}}}(s)=0\qquad\mbox{for $\alpha\in\Lambda _0$}
\end{equation}
and
\begin{equation}\label{eq77}
\rho _{_{f_{i}, g_{(\beta)1}}}(s)=0\qquad\mbox{for $\beta\in\Lambda _1$},
\end{equation}
then $s=0$.

\medskip
\underline{\it Case 1}: $i=0$, in which case, $\Lambda _0\ne \emptyset$. Since 
$s\in \mathcal{O}(D^{\sharp}(f_0))$, we have
\begin{equation}\label{eq78}
s=\left(\frac{a_0}{f_0^n}+\frac{a_1}{z_1f_0^m},\,\,
b_0+\frac{b_1}{1^{\sharp}f_0^k}\right)
\end{equation}
for some $n$, $m$, $k\in \mathcal{Z}_{\ge 0}$, $z_1f_0^m\in T_1(f_0)$, $a_i$, $b_i\in R_i$ and $i=0$, $1$. Since 
$D^{\sharp}(g_{(\alpha)0})\subseteq D^{\sharp}(f_0)$, there exist some $r_{(\alpha)0}\in R_0$ and $u_{(\alpha)}\in\mathcal{Z}_{>0}$ such that
\begin{equation}\label{eq79}
g_{(\alpha)0}^{u_{(\alpha)}}=r_{(\alpha)0}f_0\quad\mbox{for each $\alpha\in\Lambda _0$. }
\end{equation}

$\bullet$ If $D^{\sharp}(f_{0})\bigcap Spec^{\sharp}_1R\ne\emptyset$, then for each 
$\alpha\in\Lambda _0$, we have 
$$\mbox{either}\quad
D^{\sharp}(g_{(\alpha)0})\bigcap Spec^{\sharp}_1R\ne\emptyset\quad \mbox{or}\quad
D^{\sharp}(g_{(\alpha)0})\bigcap Spec^{\sharp}_1R=\emptyset.
$$
Using (\ref{eq60}), (\ref{eq592}), (\ref{eq76}) and (\ref{eq79}), we have
\begin{eqnarray}\label{eq80}
&&0=\rho _{_{f_{0}, g_{(\alpha)0}}}(s)=
\rho _{_{f_{0}, g_{(\alpha)0}}}\left(\left(\frac{a_0}{f_0^n}+\frac{a_1}{z_1f_0^m},\,\,
b_0+\frac{b_1}{1^{\sharp}f_0^k}\right)\right)\nonumber\\
&=& \left\{\begin{array}{ll}
\left(\frac{a_0r_{(\alpha)0}^n}{g_0^{nu_{(\alpha)}}}+
\frac{a_1r_{(\alpha)0}^m}{z_1g_0^{mu_{(\alpha)}}},\,\,
b_0+\frac{b_1r_{(\alpha)0}^k}{1^{\sharp}g_0^{ku_{(\alpha)}}}\right)&\mbox{if $D^{\sharp}(g_{(\alpha)0})\bigcap Spec^{\sharp}_1R\ne\emptyset$}\\&\\
\left(\frac{a_0r_{(\alpha)0}^n}{g_0^{nu_{(\alpha)}}},\,\,b_0\right)&
\mbox{if $D^{\sharp}(g_{(\alpha)0})\bigcap Spec^{\sharp}_1R=\emptyset$}
\end{array}\right.
\end{eqnarray}

$\bullet$ If $D^{\sharp}(f_{0})\bigcap Spec^{\sharp}_1R=\emptyset$, then we have 
$D^{\sharp}(g_{(\alpha)0})\bigcap Spec^{\sharp}_1R=\emptyset$ for each 
$\alpha\in\Lambda _0$. By (\ref{eq592}), (\ref{eq76}) and (\ref{eq79}), we have
\begin{eqnarray}\label{eq80a}
0&=&\rho _{_{f_{0}, g_{(\alpha)0}}}(s)=
\rho _{_{f_{0}, g_{(\alpha)0}}}\left(\left(\frac{a_0}{f_0^n}+\frac{a_1}{z_1f_0^m},\,\,
b_0+\frac{b_1}{1^{\sharp}f_0^k}\right)\right)\nonumber\\
&=& \left(\frac{a_0r_{(\alpha)0}^n}{g_0^{nu_{(\alpha)}}},\,\,b_0\right)\quad
\mbox{for each $\alpha\in\Lambda _0$}
\end{eqnarray}

Note  that
\begin{equation}\label{eq80b} 
D^{\sharp}(g_{(\alpha)0})\bigcap Spec^{\sharp}_1R=\emptyset\Longrightarrow
1^{\sharp}g_{(\alpha)0}^b=0\quad\mbox{for some $b\in \mathcal{Z}_{\ge 0}$}
\end{equation}

It follows from (\ref{eq80}), (\ref{eq80a}), (\ref{eq80b}) and Proposition~\ref{pr7.2} (i) that there exist 
$t$, $d$, $\delta _{\alpha _i}\in \mathcal{Z}_{>0}$,
$ \alpha _i\in \Lambda _0$ and $h_{(\alpha _i)0}\in R_0$  such that
\begin{equation}\label{eq81} 
f_0^t=\sum_{i=1}^d h_{(\alpha _i)0}\,\,g_{(\alpha _i)0}^{\delta _{\alpha _i}}
\end{equation}
and
\begin{equation}\label{eq82} 
f_0^ta_0=(1^{\sharp}f_0^t)\,\sharp\, a_1=b_0=(1^{\sharp}f_0^t)\,\sharp\, b_1=0.
\end{equation}

By (\ref{eq82}), $\displaystyle\frac{a_0}{f_0^n}+\displaystyle\frac{a_1}{1^{\sharp}f_0^m}=0$ and 
$b_0+\displaystyle\frac{b_1}{1^{\sharp}f_0^k}=0$. This proves $s=0$. 

\medskip
\underline{\it Case 2}: $i=1$, in which case, $s\in \mathcal{O}(D^{\sharp}(f_1))$. Thus 
$s=a_0+\displaystyle\frac{a_1}{f_1^{\sharp m}}$ for some $a_i\in R_i$ and 
$m\in\mathcal{Z}_{\ge 0}$. By $(\ref{eq61})$, $D^{\sharp}(g_{(\alpha)0)})=\emptyset $ for 
$\alpha\in\Lambda _0$ in this case. Thus, (\ref{eq75}) becomes
\begin{equation}\label{eq83} 
D^{\sharp}(f_1)=\bigcup _{\beta\in\Lambda _1}D^{\sharp}(g_{(\beta)1}).
\end{equation}
Since $D^{\sharp}(g_{(\beta)1})\subseteq D^{\sharp}(f_1)$, there exist some $r_{(\beta)1}\in R_1$ and $\theta _{(\beta)}\in\mathcal{Z}_{>0}$ such that
\begin{equation}\label{eq84}
g_{(\beta)1}^{\theta _{(\beta)}}=r_{(\beta)1}\,\sharp\,f_1\quad\mbox{for each 
$\beta\in\Lambda _1$. }
\end{equation}

Using (\ref{eq65}) and (\ref{eq84}), we have
$$
0=\rho _{_{f_1, g_{(\beta)1}}}(s)
=\rho _{_{f_1, g_{(\beta)1}}}\left(a_0+\displaystyle\frac{a_1}{f_1^{\sharp m}}\right)
=a_0+\displaystyle\frac{a_1\,\sharp\,r_{(\beta)1}^{\sharp m}}
{g_{(\beta)1}^{\sharp m\theta _{(\beta)}}},
$$
which implies that there exist 
$t'$, $d'$, $\delta _{\beta _i}'\in \mathcal{Z}_{>0}$,
$ \beta _i\in \Lambda _1$ and $h_{(\beta _i)0}\in R_0$  such that
\begin{equation}\label{eq85} 
f_1^{t'}=\sum_{i=1}^{d'} h_{(\beta _i)1}\,\,g_{(\beta _i)1}^{\delta _{\beta _i}'}
\end{equation}
and
\begin{equation}\label{eq86} 
a_0=f_1^{\sharp t'}\,\sharp\, a_1=0.
\end{equation}
Hence, $s=0$ by (\ref{eq86}).

\medskip
This proves that the properties (i) in Definition~\ref{def8.2} holds for the basic open set $D^{\sharp}(f_i)$ and the open covering (\ref{eq75}). 

\bigskip
In order to prove that the properties (ii) in Definition~\ref{def8.2} holds for the basic open set $D^{\sharp}(f_i)$ and its open covering
\begin{equation}\label{eq87}
\{\, D^{\sharp}(g_{(\alpha)0})\,|\, \alpha\in\Lambda _0\,\}\bigcup
\{\, D^{\sharp}(g_{(\beta)1})\,|\, \beta\in\Lambda _1\,\},
\end{equation}
it is sufficient to prove that the properties (ii) in Definition~\ref{def8.2} holds for the basic open set $D^{\sharp}(f_i)$ and a finite subcovering of the open covering given by (\ref{eq87}). A direct computation proves  that this fact is indeed true. 

\medskip
This completes the proof of  Proposition~\ref{pr8.2}.

\hfill\raisebox{1mm}{\framebox[2mm]{}}

\end{document}